\documentclass[11pt,reqno]{amsart}
\usepackage{amsmath,amstext,amssymb,amscd}
\usepackage{verbatim}
\usepackage{enumerate}
\usepackage{mathrsfs}
\usepackage[dvipsnames,usenames]{xcolor}

\usepackage{amsthm}

\newtheorem{theorem}{Theorem}[section]

\theoremstyle{definition}
\newtheorem{definition}[theorem]{Definition}

\theoremstyle{remark}
\newtheorem{remark}[theorem]{Remark}

\numberwithin{equation}{section}

\begin{document}

\newcommand{\sgn}{\operatorname{sgn}}

\def\a{\alpha}
\def\b{\beta}
\def\d{\delta}
\def\g{\gamma}
\def\l{\lambda}
\def\o{\omega}
\def\s{\sigma}
\def\t{\tau}
\def\th{\theta}
\def\r{\rho}
\def\D{\Delta}
\def\G{\Gamma}
\def\O{\Omega}
\def\e{\varepsilon}
\def\p{\phi}
\def\P{\Phi}
\def\S{\Psi}
\def\E{\eta}
\def\m{\mu}
\def\grad{\nabla}
\def\bar{\overline}
\newcommand{\reals}{\mathbb{R}}
\newcommand{\naturals}{\mathbb{N}}
\newcommand{\ints}{\mathbb{Z}}
\newcommand{\complex}{\mathbb{C}}
\newcommand{\rationals}{\mathbb{Q}}
\newcommand{\innerprod}[1]{\left\langle#1\right\rangle}
\newcommand{\norm}[1]{\left\|#1\right\|}
\newcommand{\abs}[1]{\left|#1\right|}

\title[Wave Equations with Supercritical Source and Damping Terms]
{Global well-posedness for nonlinear wave equations with supercritical source and damping terms}

\author{Yanqiu Guo}

\date{October 29, 2018}
\keywords{nonlinear wave equations, supercritical, source terms, damping, global well-posedness}

\maketitle

\begin{abstract}
We prove the global well-posedness of weak solutions for nonlinear wave equations with supercritical source and damping terms on a three-dimensional torus $\mathbb T^3$ of the prototype
\begin{align*}
&u_{tt}-\Delta u+|u_t|^{m-1}u_t=|u|^{p-1}u,  \;\;   (x,t)  \in \mathbb T^3  \times \mathbb R^+ ;    \notag\\
&u(0)=u_0 \in H^1(\mathbb T^3)\cap L^{m+1}(\mathbb T^3), \;\; u_t(0)=u_1\in L^2(\mathbb T^3), 
\end{align*}
where $1\leq p\leq \min\{ \frac{2}{3} m + \frac{5}{3} , m \}$. Notably, $p$ is allowed to be larger than $6$.
\end{abstract}

\section{Introduction}

\subsection{The model}
Wave equations under the influence of nonlinear damping and source terms have attracted considerable attention.
The canonical equation of this type reads
\begin{align} 
u_{tt}-\Delta u+  |u_t|^{m-1} u_t  =  |u|^{p-1} u,     \label{mo}
\end{align}
where $m, p\geq 1$. A major interest for this topic lies in understanding the ``competition" between the frictional damping term $|u_t|^{m-1} u_t$ and the energy-amplifying source term $|u|^{p-1} u$.

The purpose of this paper is to provide a suitable assumption on $p$ and $m$, such that model (\ref{mo}) is globally well-posed for weak solutions defined in a three-dimensional periodic physical domain, and the source term $|u|^{p-1}u$ is allowed to have a ``fast" growth rate $p\geq 6$.

Let us review some important results in the literature which are concerned with equation (\ref{mo}).
Georgiev and Todorova \cite{Geo-Todo} studied (\ref{mo}) in a bounded domain $\Omega \subset \mathbb R^3$ with a Dirichl\'et boundary condition. For a source term of subcritical or critical growth rate ($1\leq p\leq 3$), they proved the well-posedness of weak solutions for (\ref{mo}). In addition, the solution is global if the damping dominates the source term in the sense that $m\geq p$, whereas the solution blows up in finite time if the strength of the source exceeds the damping effect, namely, $p>m$.

We remark that $p=3$ is called the critical growth rate for the source term $|u|^{p-1}u$ because the operator $u \mapsto u^3$ is locally Lipschitz continuous from $H^1$ to $L^2$ in three dimensions.

Bociu and Lasiecka \cite{BL1,BL2, BL3} considered (\ref{mo}) with supercritical source terms, in a bounded domain $\Omega \subset \mathbb R^3$ satisfying a Newman boundary condition, and showed the existence and uniqueness of weak solutions if $1\leq p\leq \frac{6m}{m+1}$, allowing $p$ have the range $[1,6)$.

In the literature, the Cauchy problem for (\ref{mo}) in $\mathbb R^n$ was also investigated (see, e.g., \cite{Radu, To}). Moreover, it is of interest to consider interaction between source terms and other types of damping terms in nonlinear wave equations, for instance, strong damping (e.g., \cite{GS}), degenerate damping (e.g., \cite{BLR,BLR2}), and viscoelastic damping (e.g. \cite{GRSTT, GRS, GRS1}). One may also refer to \cite{AR, BGRT, CCL, GR, GR1, GR2, LS, STV, Viti, Viti2} and references therein for more works on nonlinear wave equations with damping and source terms. It is also worth mentioning papers \cite{YT, KN} on analyticity for a class of nonlinear wave equations including (\ref{mo}) as a special case.

\vspace{0.1 in}

\subsection{Main results}
In this paper, we study the following nonlinear wave equation with damping and source terms 
defined in a three-dimensional fundamental periodic domain $\mathbb T^3 = [-\pi,\pi]^3$:
\begin{align} 
&u_{tt}-\Delta u+  |u_t|^{m-1} u_t  =  |u|^{p-1} u,   \;\; (x,t)  \in \mathbb T^3  \times \mathbb R^+;     \label{sys1}\\    
&u(x,0)=u_0(x)\in  H^1(\mathbb T^3)\cap L^{m+1}(\mathbb T^3),    \;\;  u_t(x,0)=u_1(x) \in L^2(\mathbb T^3),            \label{sys2}
\end{align}
where $m, p\geq 1$.    Our main result states, 
if $1\leq p\leq \min\{ \frac{2}{3} m + \frac{5}{3} , m\}$, then system (\ref{sys1})-(\ref{sys2}) admits a unique global weak solution which depends continuously on initial data. Note, in the initial condition (\ref{sys2}), we demand an extra integrability for $u_0$, namely, $u_0 \in L^{m+1}(\mathbb T^3)$ if $m>5$.

We choose the physical domain to be a torus $\mathbb T^3$ because we want to focus on the interaction between the damping and source terms, without influence of boundary conditions. Also, we restrict our analysis to 3D since it is more physically relevant. Our results extend easily to an $n$-dimensional torus $\mathbb T^n$, by accounting for the corresponding Sobolev imbeddings, and accordingly adjusting the conditions imposed on the parameters.

Throughout the paper, we denote by $\|\cdot \|_s =  \|\cdot\|_{L^s(\mathbb T^3)}$ for $L^s$-norm. Also, for a function $y(x,t)$ defined on $\mathbb T^3 \times \mathbb R^+$, the partial derivative in $t$ is denoted by $y'=y_t=\frac{\partial y}{\partial t}$.

Let us introduce the definition of a weak solution for system (\ref{sys1})-(\ref{sys2}).

\begin{definition}   \label{weak}
Let $T>0$. We call $(u,u_t)$ a \emph{weak solution} for system (\ref{sys1})-(\ref{sys2}) on $[0,T]$ if 
\begin{itemize}
\item   $u(x,0)=u_0(x) \in H^1(\mathbb T^3)  \cap L^{m+1}(\mathbb T^3)$, $u_t(x,0)=u_1(x)\in L^2(\mathbb T^3)$;
\item $u \in L^{\infty}(0,T;H^1(\mathbb T^3)) \cap  L^{\infty}(0,T;L^{m+1}(\mathbb T^3)) $; \\
$u_t \in L^{\infty}(0,T;L^2(\mathbb T^3)) \cap L^{m+1}(\mathbb T^3 \times (0,T))$;  \\
 $u_{tt} \in L^{\frac{m+1}{m}}(0,T; X')$ where $X= H^1(\mathbb T^3) \cap L^{m+1}(\mathbb T^3)$;
\item $(u,u_t)$ verifies the identity
\begin{align}  \label{varI}   
&\int_{\mathbb T^3} u_t(t) \phi(t) dx - \int_{\mathbb T^3} u_t(0) \phi(0) dx - \int_0^t \int_{\mathbb T^3} u_t(\tau) \phi_t (\tau) dx d\tau 
+  \int_0^t  \int_{\mathbb T^3} \nabla u (\tau)  \cdot \nabla  \phi (\tau) dx d\tau  \notag\\
& +   \int_0^t \int_{\mathbb T^3}   |u_t(\tau)|^{m-1} u_t(\tau)   \phi(\tau) dx d\tau = \int_0^t   \int_{\mathbb T^3} |u(\tau)|^{p-1} u(\tau) \phi(\tau) dx d\tau,
\end{align}
for all $t\in [0,T]$, and for any $\phi \in C([0,T]; H^1(\mathbb T^3)) \cap L^{m+1}(\mathbb T^3 \times (0,T))$ with $\phi_t \in C([0,T];L^2(\mathbb T^3))$.
\end{itemize}
\end{definition}

Our first theorem deals with the global existence of weak solution for the initial value problem (\ref{sys1})-(\ref{sys2}). Also the energy identity  holds for weak solutions. Moreover, a global solution $(u,u_t)$ grows at most exponentially in time. 

\begin{theorem}[Global existence of weak solutions]  \label{thm1} 
Assume either \emph{Case 1}: $1\leq p\leq m \leq 5$ or \emph{Case 2}: $1\leq p<\frac{5}{6}(m+1)$ for $m>5$.
Suppose $u_0 \in H^1(\mathbb T^3)\cap L^{m+1}(\mathbb T^3)$ and $u_1 \in L^2(\mathbb T^3)$.
Let $T>0$ be arbitrarily large. Then, system (\ref{sys1})-(\ref{sys2}) has a weak solution $(u,u_t)$ on $[0,T]$ in the sense of Definition \ref{weak}.
Also, the energy identity holds:
\begin{align}  \label{EI}
E(t) +  \int_0^t \|u_t(\tau)\|^{m+1}_{m+1}d\tau = E(0), \;\;  \text{for all} \;\;   t\in [0,T],
\end{align}
where the total energy $E(t):=     \frac{1}{2} (\|\nabla u(t)\|_2^2 + \|u_t(t)\|_2^2) - \frac{1}{p+1} \|u(t)\|_{p+1}^{p+1} $. In addition, 
\begin{align}   \label{expm}
\mathscr E(t)  +   \frac{1}{2} \int_0^t  \|u_t(\tau)\|^{m+1}_{m+1} d\t \leq  \left(\mathscr E(0) + t\right)      e^{Ct}, \;\;   \text{for all}  \;\;  t\in [0,T],
\end{align}
where $\mathscr E(t):=    \frac{1}{2}\left(\|\nabla u(t)\|^2_2+\|u_t(t)\|^2_2 \right)+\frac{1}{m+1}\norm{u(t)}_{m+1}^{m+1}$.
\end{theorem}

Our second theorem establishes the uniqueness of weak solutions by assuming a slightly stronger restriction on $(m,p)$. Continuous dependence on initial data is also provided.
\begin{theorem}[Uniqueness and continuous dependence]  \label{thm2}
Assume either \emph{Case I}: $1\leq p\leq m \leq 5$ or \emph{Case II}: $1\leq p \leq \frac{2}{3} m + \frac{5}{3} $ for $m>5$.
Suppose $u_0 \in H^1(\mathbb T^3)\cap L^{m+1}(\mathbb T^3)$ and $u_1 \in L^2(\mathbb T^3)$. Let $T>0$ be arbitrarily large. Then, system (\ref{sys1})-(\ref{sys2}) has a unique weak solution $(u,u_t)$ on $[0,T]$ in the sense of Definition \ref{weak}. 
Also, the weak solution depends continuously on initial data. More precisely, let $(u_0^n,u_1^n)$ be a sequence of initial data such that
$\lim_{n \rightarrow \infty}\|u_0^n - u_0\|_{H^1}=0$, $\lim_{n \rightarrow \infty}\|u_0^n - u_0\|_{m+1}=0$ and $\lim_{n \rightarrow \infty}\|u_1^n - u_1\|_2=0$, then the corresponding sequence of weak solutions $(u_n,u_n')$ converges to $(u,u_t)$ in the sense that
\begin{align}   \label{cdconverge}
\lim_{n\rightarrow \infty} \left[ \sup_{t\in [0,T]} \left( \|u_n-u\|_{H^1}^2     + \|u_n-u\|_{m+1}^{m+1}     +   \|u_n'- u_t\|_2^2    \right)\right] =0.
\end{align}
\end{theorem}
\begin{remark}    \label{remark2}
The range of $(m,p)$ assumed in Theorem \ref{thm2} (i.e. the ``union" of Case I and Case II) can be equivalently expressed as $1\leq p\leq \min\{ \frac{2}{3} m + \frac{5}{3} , m \}$.
Also, Case II of Theorem \ref{thm2} is a slightly smaller region in the $(m,p)$ plane compared to Case 2 of Theorem \ref{thm1}.
\end{remark}

\begin{remark}     
For the sake of clarity, we consider the ``typical" frictional damping term $|u_t|^{m-1}u_t$ and the ``typical" source term $|u|^{p-1}u$. Nonetheless, our results hold for more general damping and source terms. More precisely, the damping term  $|u_t|^{m-1}u_t$ can be generalized to $g(u_t)$ where $g\in C(\mathbb R)$ is a monotone increasing function vanishing at the origin such that 
$$a|s|^{m+1} \leq g(s)s\leq b|s|^{m+1},   \;\; \text{where}  \;\;  b\geq a>0   \;\; \text{and}   \;\;  m\geq 1.$$
Also, the source term $|u|^{p-1} u$ can be generalized to $h(u)$ where $h$ is a $C^1(\mathbb R)$ function ($C^2$ is required if $p>3$) satisfying
 \begin{align*}
 \begin{cases}
 |h'(s)|\leq C(|s|^{p-1} +1),  \;\; \text{if} \;\;1\leq p \leq 3;\\
 |h''(s)| \leq C(|s|^{p-2}  +1),   \;\; \text{if}  \;\; p>3.
 \end{cases}
\end{align*}
\end{remark}

\vspace{0.1 in}

\section{Global existence of weak solutions}
This section is devoted to proving Theorem \ref{thm1}, namely, the existence of global weak solutions, the energy identity, and the exponential bound for the growth of solutions at large time.  

\subsection{Galerkin approximation system}
We show the existence of weak solutions for system (\ref{sys1})-(\ref{sys2}) via the standard Galerkin approximation method.
Let us first review some classical results regarding Fourier series on a torus. 
For a periodic function $f\in L^1(\mathbb T^3)$ where $\mathbb T^3=[-\pi,\pi]^3$, the $k$th Fourier coefficient of $f$ is defined by
$\hat f(k) =  \int_{\mathbb T^3}  f(x) e^{- i k\cdot x }dx$. The Fourier series of $f$ at $x\in \mathbb T^3$ is written as $\sum_{k\in \mathbb Z^3} \hat f(k) e^{i k \cdot x}$. We define the square partial sum of the Fourier series of $f$ by
\begin{align}  \label{Spartial}
P_n f(x)= \sum_{\substack{k=(k_1,k_2,k_3) \in \mathbb Z^3 \\ |k_1|, |k_2|, |k_3| \leq n}} \hat f (k) e^{ik \cdot x}.
\end{align}
Note that, for a Fourier series, the square partial sum defined in (\ref{Spartial}) is contrast to the spherical partial sum:
 $\sum_{|k| \leq n} \hat f (k) e^{ik \cdot x}$. 
It is a classical result that the square partial sum $P_n f$ converges to $f$ in $L^s(\mathbb T^3)$ for any $s\in (1,\infty)$ 
(see, e.g., \cite{Grafakos}), namely, for an $f\in L^s(\mathbb T^3)$ with $1<s<\infty$,
\begin{align}  \label{PSC}
\lim_{n\rightarrow \infty} \|P_n f  - f\|_{L^{s}(\mathbb T^3)} =0.
\end{align}
Moreover, for any $f\in L^s(\mathbb T^3)$ with $1<s<\infty$, 
\begin{align}   \label{PSC1}
\|P_n f\|_{L^s(\mathbb T^3)} \leq c_s \|f\|_{L^s(\mathbb T^3)},
\end{align}
for some positive constant $c_s$ independent of $n$ and $f$.

We consider the Galerkin approximation system
\begin{align}  
&u_n''-\Delta u_n+P_n (|u_n'|^{m-1} u_n')=P_n (|u_n|^{p-1}u_n),    \;\;  (x,t) \in  \mathbb T^3 \times \mathbb R^+;  \label{Gal}    \\ 
&u_n(x,0)=P_n u_0(x) , \;\; u_n'(x,0) = P_n u_1(x),   \label{Gal-1}
\end{align}
where $u_n(x,t)$ is a trigonometric polynomial of the form:
$$u_n(x,t)=\sum_{\substack{k=(k_1,k_2,k_3) \in \mathbb Z^3\\ |k_1|, |k_2|, |k_3| \leq n}} \hat u_n(k,t) e^{ik\cdot x}.$$

By the Cauchy-Peano theorem, for each $n$, Galerkin system (\ref{Gal})-(\ref{Gal-1}) has a solution $u_n$ on $[0,T_n)$ for some $T_n\in (0,\infty]$ which stands for the maximum life span.

\vspace{0.1 in}

\subsection{Energy estimate}
In this subsection, we show that $(u_n,u_n')$ is bounded in the energy space $H^1(\mathbb T^3) \times L^2(\mathbb T^3)$ uniformly in $n$.
Multiply (\ref{Gal}) by $u_n'$ and integrate over $\mathbb T^3 \times (0,t)$. One has, for $t\in [0,T_n)$,
\begin{align}   \label{ener-idn0}
&\frac{1}{2} \left( \|\nabla u_n(t)\|_2^2   +  \|u_n'(t)\|_2^2       \right)  +       \int_0^t \int_{\mathbb T^3}  |u_n'|^{m+1}  dx d\tau      \notag\\
&=   \frac{1}{2} \left(   \|\nabla u_n(0)\|_2^2   + \|u_n'(0)\|_2^2     \right)   +    \int_0^t \int_{\mathbb T^3} |u_n|^{p-1} u_n  u_n' dx d\tau.
\end{align}

Define a \emph{modified} energy:
\begin{align}  \label{energy-mo}
\mathscr E_n(t)=\frac{1}{2}\left(\|\nabla u_n(t)\|^2_2     +      \|u_n'(t)\|^2_2        \right)+\frac{1}{m+1}   \norm{u_n(t)}_{m+1}^{m+1} \geq 0.
\end{align}

Then (\ref{ener-idn0}) can be written as
\begin{align}   \label{ener-idn1}
&  \mathscr E_n(t)    +       \int_0^t \int_{\mathbb T^3} |u_n'|^{m+1}   dx d\tau      \notag\\
&=  \mathscr E_n(0)   +    \int_0^t \int_{\mathbb T^3}  |u_n|^{p-1}  u_n  u_n' dx d\tau  +    \frac{1}{m+1}  \int_{\mathbb T^3} \left( |u_n(t)|^{m+1}  -   |u_n(0)|^{m+1}\right) dx.     \end{align}
Since, for $m\geq 1$, 
\begin{align*} 
&\frac{1}{m+1}  \int_{\mathbb T^3} \left( |u_n(t)|^{m+1}  -   |u_n(0)|^{m+1}\right) dx  \notag\\
&=\int_{\mathbb T^3} \int_0^t     \frac{d}{d\tau} \left( \frac{1}{m+1} |u_n(\tau)|^{m+1} \right)   d\tau dx  
=   \int_{\mathbb T^3} \int_0^t    |u_n(\tau)|^{m-1} u_n(\tau) u_n'(\tau)   d\tau dx,
\end{align*}  
it follows from (\ref{ener-idn1}) that 
\begin{align} \label{ener-idn}
&\mathscr E_n(t) +\int_0^t \int_{\mathbb T^3} |u_n'|^{m+1} dx d\t  \notag\\
&=\mathscr E_n(0)+\int_0^t \int_{\mathbb T^3} |u_n|^{p-1} u_n u_n' dx d\t+\int_0^t \int_{\mathbb T^3}|u_n|^{m-1}u_n u_n' dx d\t.
\end{align}

We estimate the integrals on the right-hand side of (\ref{ener-idn}). Since $m\geq p\geq 1$, by using H\"older's inequality and Young's inequality, we obtain
\begin{align}   \label{fast-1}
&\int_0^t \int_{\mathbb T^3} |u_n|^p |u_n'| dx d\t \leq   Ct^{\frac{m-p}{m+1}}
\left(\int_0^t \int_{\mathbb T^3}  |u_n|^{m+1}  dx d\tau \right)^{\frac{p}{m+1}} \left(\int_0^t \int_{\mathbb T^3}  |u_n'|^{m+1}  dx d\tau \right)^{\frac{1}{m+1}}   \notag\\ 
&\leq \frac{1}{4}\int_0^t \int_{\mathbb T^3}  |u_n'|^{m+1}  dx d\tau
+  C  \int_0^t \int_{\mathbb T^3}  |u_n|^{m+1}  dx d\tau +   t.
\end{align}
Similarly, we have 
\begin{align}   \label{fast-2}
\int_0^t \int_{\mathbb T^3}  |u_n|^m  |u_n'|  dx d\tau \leq  \frac{1}{4}   \int_0^t  \int_{\mathbb T^3} |u_n'|^{m+1} dx d\tau   +  C  \int_0^t \int_{\mathbb T^3} |u_n|^{m+1} dx d\tau.
 \end{align}
Substituting (\ref{fast-1})-(\ref{fast-2}) into (\ref{ener-idn}) yields
 \begin{align*}
\mathscr E_n(t)  +   \frac{1}{2} \int_0^t \int_{\mathbb T^3} |u_n'|^{m+1} dx d\t  
&\leq \mathscr E_n(0) +    C \int_0^t   \int_{\mathbb T^3} |u_n|^{m+1} dx d\tau  +  t  \notag\\
&\leq    \mathscr E_n(0) +    C \int_0^t \mathscr E_n(\tau) d\tau  +  t,     \;\;\text{for all} \,\, t\in [0,T_n).
\end{align*}
Then, using Gr\"onwall's inequality, we obtain
\begin{align}   \label{exp-grow}
\mathscr E_n(t)  +   \frac{1}{2} \int_0^t \int_{\mathbb T^3} |u_n'|^{m+1} dx d\t
\leq  \left(\mathscr E_n(0) + t\right)      e^{Ct},    \;\; \text{for all}  \;\; t\in [0,T_n).
\end{align}

By (\ref{energy-mo}), $$\mathscr E_n(0)= \frac{1}{2} (\|\nabla (P_n u_0)\|_2^2  +  \| P_n u_1\|_2^2) +  \frac{1}{m+1} \|P_n u_0\|_{m+1}^{m+1}.$$
According to Plancherel's theorem, we have $\| P_n u_1\|_2^2 \leq \|u_1\|_2^2 $ for all $n\in \mathbb N$  
and $ \| \nabla (P_n u_0)\|_2^2  =   \| P_n (\nabla u_0)\|_2^2 \leq \|\nabla u_0\|_2^2$ for all $n\in \mathbb N$.
Moreover, since $u_0\in L^{m+1}(\mathbb T^3)$ with $m\geq 1$, by virtue of (\ref{PSC1}), 
we have $\|P_n u_0\|_{m+1} \leq  c_m \|u_0 \|_{m+1}$, for some constant $c_m$ independent of $n$ and $u_0$.
Therefore, $\mathscr E_n(0)$ has an upper bound independent of $n$, namely
\begin{align} \label{sEn0}
\mathscr E_n(0) \leq  \frac{1}{2} (\|\nabla u_0\|_2^2 + \|u_1\|_2^2) + c_m \|u_0\|_{m+1}^{m+1}.
\end{align}

Applying (\ref{sEn0}) to the right-hand side of (\ref{exp-grow}) yields, for any $n\in \mathbb N$,
\begin{align}   \label{exp-growth}
&\frac{1}{2}\left(\|\nabla u_n(t)\|^2_2+\|u_n'(t)\|^2_2 \right)+\frac{1}{m+1}\norm{u_n(t)}_{m+1}^{m+1} 
 +   \frac{1}{2} \int_0^t   \|u_n'(\tau)\|^{m+1}_{m+1} d\t    \notag\\
&\leq    \left(\frac{1}{2} (\|\nabla u_0\|_2^2 + \|u_1\|_2^2) + c_m \|u_0\|_{m+1}^{m+1} + t \right)e^{Ct},
\end{align}
for all $t\in [0,T_n)$. Since $u_0 \in H^1(\mathbb T^3) \cap L^{m+1}(\mathbb T^3)$ and $u_1\in L^2(\mathbb T^3)$, we
see from (\ref{exp-growth}) that a solution $(u_n,u_n')$ for the Galerkin system (\ref{Gal})-(\ref{Gal-1}) does not blow up at $T_n$. Therefore, we conclude $T_n=\infty$, namely, 
$(u_n,u_n')$ is a global solution for (\ref{Gal}) for all $t\in [0,\infty)$.

For $m\geq 1$, we have
\begin{align}  \label{pH1}
\|u_n\|_{H^1(\mathbb T^3)}^2  =  \|\nabla u_n\|_2^2 + \|u_n\|_2^2  
\leq  \|\nabla u_n\|_2^2 +  C\|u_n\|_{m+1}^{m+1} + 1.
\end{align}
Let $T>0$. Since (\ref{exp-growth}) holds for all $t\geq 0$, then by (\ref{pH1}), one has
\begin{align}     \label{ub-1}
u_n   \;\; \text{is uniformly bounded in}  \;\;  L^{\infty}(0,T; H^1(\mathbb T^3)).
\end{align}
Moreover, by (\ref{exp-growth}), we have
\begin{align}
&u_n   \;\; \text{is uniformly bounded in}  \;\;  L^{\infty}(0,T; L^{m+1}(\mathbb T^3));    \label{ub-2} \\
&u_n'   \;\; \text{is uniformly bounded in}  \;\;  L^{\infty}(0,T; L^2(\mathbb T^3));  \label{ub-3} \\
&u_n'   \;\; \text{is uniformly bounded in}  \;\;  L^{m+1}(\mathbb T^3 \times (0,T)).    \label{ub-4}
\end{align}

Notice $\int_{\mathbb T^3} \left|   |u_n|^{p-1} u_n     \right|^{\frac{m+1}{m}} dx 
=\int_{\mathbb T^3}   |u_n|^{\frac{(m+1)p}{m}}  dx   \leq    C \left(\int_{\mathbb T^3}  |u_n|^{m+1} dx\right)^{\frac{p}{m}}$ for $m\geq p \geq 1$. Thus, because of (\ref{ub-2}), we obtain
\begin{align}   \label{unif-bd}
 |u_n|^{p-1} u_n  \;\;  \text{is uniformly bounded in}\;\;  L^{\infty}(0,T;L^{\frac{m+1}{m}}(\mathbb T^3)).
 \end{align}

Also, since $\int_0^T \int_{\mathbb T^3} \left| |u_n'|^{m-1} u_n' \right|^{\frac{m+1}{m}} dx dt  
=     \int_0^T \int_{\mathbb T^3} |u_n'|^{m+1} dx dt $, and due to (\ref{ub-4}),
it follows that 
\begin{align}    \label{unif-bdt}
|u_n'|^{m-1} u_n'  \;\;  \text{is uniformly bounded in}   \;\; L^{\frac{m+1}{m}}(\mathbb T^3 \times (0,T)).
\end{align}

Because of the Galerkin equation (\ref{Gal}), we have $u_n''=\Delta u_n - P_n (|u_n'|^{m-1} u_n')+P_n (|u_n|^{p-1}u_n)$. Note that
$\Delta u_n$ is uniformly bounded in $L^{\infty}(0,T; (H^1(\mathbb T^3))')$ due to (\ref{ub-1}). Thus, by virtue of (\ref{PSC1}), (\ref{unif-bd}) and (\ref{unif-bdt}), we obtain 
\begin{align}  \label{ub-5}
u_n''    \;\;   \text{is uniformly bounded in} \;\;  L^{\frac{m+1}{m}}(0,T; X'),  \;\; \text{where} \;\;  X=H^1(\mathbb T^3) \cap L^{m+1}(\mathbb T^3).
\end{align}

By virtue of the uniform bounds (\ref{ub-1})-(\ref{ub-4}) and (\ref{ub-5}), there exists a subsequence of $u_n$ satisfying 
\begin{align}
&u_n \rightarrow u   \;\;  \text{weakly$^*$ in} \;\;   L^{\infty} (0,T; H^1(\mathbb T^3));   \label{con-1}  \\
&u_n \rightarrow u         \;\;  \text{weakly$^*$ in} \;\;     L^{\infty}(0,T;   L^{m+1}(\mathbb T^3));   \label{con-2}  \\
&u_n' \rightarrow u'    \;\;  \text{weakly$^*$ in} \;\;    L^{\infty} (0,T; L^2(\mathbb T^3));      \label{con-3}  \\
&u_n' \rightarrow u'    \;\;  \text{weakly in} \;\;    L^{m+1} (\mathbb T^3 \times (0,T));   \label{con-4}  \\
&u_n'' \rightarrow u''     \;\;  \text{weakly$^*$ in} \;\;     L^{\frac{m+1}{m}}(0,T; X'),     \label{con-5}
\end{align}
where  $X=H^1(\mathbb T^3) \cap L^{m+1}(\mathbb T^3)$. 

Moreover, because of (\ref{ub-1}) and (\ref{ub-3}), and by the compact imbedding $H^1 \hookrightarrow H^{1-\epsilon} \hookrightarrow L^2$ for an $\epsilon\in (0,1)$, 
we conclude from the Aubin-Lions-Simon lemma that, on a subsequence
\begin{align}    \label{strong}
u_n \rightarrow u  \;\;   \text{in}  \;\;   C([0,T]; H^{1-\epsilon}(\mathbb T^3)).
\end{align}

Furthermore, due to (\ref{strong}), one can extract a subsequence
\begin{align}    \label{ptwise}
u_n  \rightarrow  u    \;\;  \text{almost everywhere in} \;\;  \mathbb T^3 \times (0,T).
\end{align}

\vspace{0.1 in}

\subsection{Convergence of the source term $|u_n |^{p-1}u_n$} \label{sec-cs}
We show $|u_n |^{p-1}u_n$ converges weakly to $|u |^{p-1}u$ in $L^{\frac{m+1}{m}}(\mathbb T^3 \times (0,T))$.

In fact,  $\left| |u_n|^{p-1} u_n   -      |u|^{p-1} u  \right|   \leq C  (|u_n|^{p-1}   +   |u|^{p-1})  |u_n - u|$, for $p\geq 1$. Thus, due to (\ref{ptwise}), we have
\begin{align}    \label{faec}
|u_n|^{p-1} u_n   \rightarrow  |u|^{p-1} u \;\;  \text{almost everywhere in} \;\;  \mathbb T^3 \times (0,T).
\end{align}

Recall a real analysis result: for a sequence of functions $f_n $ defined on a measure space $Y$, if $\sup_{n}   \|f_n\|_{L^s(Y)} <\infty$
and $f_n \rightarrow f$ a.e. in $Y$, then  $f_n \rightarrow f$ weakly in $L^s(Y)$ if $1<s<\infty$ (see, e.g., \cite{Folland}).
Here, by (\ref{unif-bd}), we know $|u_n|^{p-1} u_n$ is uniformly bounded in $L^{\frac{m+1}{m}}(\mathbb T^3 \times (0,T))$, thus along with (\ref{faec}), we conclude
\begin{align}  \label{conv-fn}
 |u_n|^{p-1} u_n \rightarrow |u|^{p-1} u  \text{\;\;weakly in \;\;} L^{\frac{m+1}{m}}(\mathbb T^3 \times (0,T)).
\end{align}

\vspace{0.1 in}

\subsection{Convergence of the damping term $|u_n'|^{m-1} u_n'$}
In this section, we show that $|u_n'|^{m-1} u_n'$ converges weakly to $|u'|^{m-1} u'$ in $L^{\frac{m+1}{m}}(\mathbb T^3 \times (0,T))$. 
The monotonicity of the damping term is critical to our argument.

Thanks to (\ref{unif-bdt}), there exists a function $\psi \in L^{\frac{m+1}{m}}(\mathbb T^3 \times (0,T))$ and a subsequence of $|u_n'|^{m-1} u_n'$ such that
\begin{align} \label{g-conv}
|u_n'|^{m-1} u_n' \longrightarrow \psi \text{\;\;weakly in\;\;} L^{\frac{m+1}{m}}(\mathbb T^3 \times (0,T)).
\end{align}
It remains to show $\psi  =   |u'|^{m-1} u' $.

Set $w=u_n-u_j$. Due to the Galerkin system (\ref{Gal}), the following equality is valid.
\begin{align} \label{diff-ener0}
&\frac{1}{2}\left(  \|\nabla w(t)\|^2_2    +     \|w'(t)\|^2_2      \right)+\int_0^t \int_{\mathbb T^3} \left[P_n(|u_n'|^{m-1} u_n')-P_j(|u_j'|^{m-1} u_j')\right]  w' dx d\t   \notag\\
&=\frac{1}{2}\left(\|\nabla w(0)\|^2_2   +      \|w'(0)\|^2_2          \right)+ \int_0^t \int_{\mathbb T^3} \left[ P_n(|u_n|^{p-1} u_n )-P_j ( |u_j|^{p-1}  u_j  )   \right] w' dx d\t.
\end{align}

We remark that the projection $P_n$ in the Galerkin system affects the monotonicity of the nonlinear damping term. 
Especially, $\int_{\mathbb T^3} \left[P_n(|u_n'|^{m-1} u_n')-P_j(|u_j'|^{m-1} u_j')\right]  w' dx$ is not necessarily positive.
To remedy the situation, we split this integral into a positive part and a ``residue" part. 
More precisely, by assming $n\geq j$, we have
\begin{align}     \label{posve}
&\int_{\mathbb T^3} \left[P_n(|u_n'|^{m-1} u_n')-P_j(|u_j'|^{m-1} u_j')\right]  w' dx 
=  \int_{\mathbb T^3} \left[|u_n'|^{m-1} u_n'  -P_j(|u_j'|^{m-1} u_j')\right]  w' dx   \notag\\
&=   \int_{\mathbb T^3} \left(   |u_n'|^{m-1} u_n'  - |u_j'|^{m-1} u_j' \right)  w' dx
+   \int_{\mathbb T^3} \left[(P_n-P_j)   (|u_j'|^{m-1} u_j' )  \right]   w' dx.
\end{align}
In addition, for the sake of convenience, we also split the integral of source terms in (\ref{diff-ener0}) in the same manner as (\ref{posve}).
As a result, if $n\geq j$, then equality (\ref{diff-ener0}) can be written as
\begin{align}   \label{diff-ener}
&\frac{1}{2}\left(   \|\nabla w(t)\|^2_2 + \|w'(t)\|^2_2 \right) +\int_0^t \int_{\mathbb T^3} \left(|u_n'|^{m-1} u_n'- |u_j'|^{m-1} u_j'\right)   w' dx d\t  \notag\\
&  \hspace{2 in}  +\int_0^t \int_{\mathbb T^3} \left[(P_n-P_j)   (|u_j'|^{m-1} u_j' )  \right]   w' dx d\t \notag\\
& = \frac{1}{2}\left(  \|\nabla w(0)\|^2_2     +     \|w'(0)\|^2_2      \right)+
\int_0^t \int_{\mathbb T^3} \left(|u_n|^{p-1} u_n - |u_j|^{p-1}  u_j \right)   w' dx d\t  \notag\\
&   \hspace{2 in}     +\int_0^t \int_{\mathbb T^3}  \left[(P_n -P_j)    (|u_j|^{p-1}  u_j)   \right]  w' dx d\t.
\end{align}
Note $(|u_n'|^{m-1} u_n'- |u_j'|^{m-1} u_j')   w' =    (|u_n'|^{m-1} u_n'- |u_j'|^{m-1} u_j')   (u_n' - u_j') \geq 0$ due to monotonicity of the function $|s|^{m-1}s$ for $m\geq 1$. Then we obtain from (\ref{diff-ener}) that, for $n\geq j$,
\begin{align}   \label{diff-ener1}
0&\leq      \frac{1}{2}\left(   \|\nabla w(t)\|^2_2 + \|w'(t)\|^2_2 \right)  +        \int_0^t \int_{\mathbb T^3} \left(|u_n'|^{m-1} u_n'- |u_j'|^{m-1} u_j'\right)   w' dx d\t  \notag\\
& \leq \frac{1}{2}\left(  \|\nabla w(0)\|^2_2   +    \|w'(0)\|^2_2     \right)+
\left|\int_0^t \int_{\mathbb T^3} \left(|u_n|^{p-1} u_n - |u_j|^{p-1}  u_j \right)   w' dx d\tau \right|  \notag\\
&  + \left|\int_0^t \int_{\mathbb T^3}  \left[(P_n -P_j)    (|u_j|^{p-1}  u_j)   \right]  w' dx d\tau \right|
+  \left|\int_0^t \int_{\mathbb T^3} \left[(P_n-P_j)   (|u_j'|^{m-1} u_j' )  \right]   w' dx d\tau\right|.
\end{align}

\subsubsection{Estimate for the ``residue" terms}

We estimate the two ``residue" terms in (\ref{diff-ener1}). They are
$\int_0^t \int_{\mathbb T^3} [(P_n-P_j)   (|u_j'|^{m-1} u_j' ) ]   w' dx d\t $ and   
$\int_0^t \int_{\mathbb T^3}  \left[(P_n -P_j)    (|u_j|^{p-1}  u_j)   \right]  w' dx d\t$, for $n\geq j$. We aim to show that they approach zero when $n$ and $j$ are large. The estimates for these two integrals are essentially the same. Thus we present the estimate for $\int_0^t \int_{\mathbb T^3} [(P_n-P_j)   (|u_j'|^{m-1} u_j' ) ]   w' dx d\t $ in details only.

Since we assume $n\geq j$ in (\ref{diff-ener}), then $P_n u'_j=u'_j=P_j u'_j$, i.e., $(P_n-P_j)u'_j=0$. 
Also, recall $w=u_n - u_j$. As a result,
\begin{align}  \label{diff-damp}
&\int_0^t \int_{\mathbb T^3} \left[(P_n-P_j)      (|u_j'|^{m-1} u_j' )  \right] w' dx d\t  
=\int_0^t \int_{\mathbb T^3}    (|u_j'|^{m-1} u_j')  \left[(P_n-P_j) (u_n'-u_j')\right] dx d\tau   \notag\\
&=\int_0^t \int_{\mathbb T^3}  (|u_j'|^{m-1} u_j')  \left[(P_n-P_j) u_n' \right] dx d\tau    \notag\\
&=\int_0^t \int_{\mathbb T^3} (|u_j'|^{m-1} u_j')  u_n'  dx d\t  -  \int_0^t \int_{\mathbb T^3} (|u_j'|^{m-1} u_j')  P_j u_n'  dx d\tau.
\end{align}
Recall that $|u_j'|^{m-1} u_j'\rightarrow \psi$ weakly in $L^{\frac{m+1}{m}}(\mathbb T^3 \times (0,T))$ by (\ref{g-conv}), 
and $u'_n \rightarrow u'$ weakly in $L^{m+1}(\mathbb T^3 \times (0,T))$ by (\ref{con-4}). Hence, for $0\leq t \leq T$,
\begin{align}   \label{ppre-bre}
\lim_{j\rightarrow \infty}   \lim_{n\rightarrow \infty} \int_0^t \int_{\mathbb T^3} (|u_j'|^{m-1} u_j') u_n'  dx d\tau  
&=   \lim_{j\rightarrow \infty}  \int_0^t \int_{\mathbb T^3} (|u_j'|^{m-1} u_j') u'  dx d\tau    \notag\\
&= \int_0^t \int_{\mathbb T^3}  \psi u'  dx d\tau.
\end{align}

Next, we look at the second term on the right-hand side of (\ref{diff-damp}). Owing to (\ref{PSC1}), we have
$\int_0^T   \int_{\mathbb T^3} |P_j (|u_j'|^{m-1} u_j' )|^{\frac{m+1}{m}} dx dt  \leq  c_m  \int_0^T   \int_{\mathbb T^3} \big| |u_j'|^{m-1} u_j' \big|^{\frac{m+1}{m}} dx dt 
=  c_m  \int_0^T  \int_{\mathbb T^3} |u_j'|^{m+1} dx dt  < \infty$ due to (\ref{ub-4}).  Hence $P_j (|u_j'|^{m-1} u_j') \in L^{\frac{m+1}{m}} (\mathbb T^3 \times (0,T))$.
Consequently, since $u'_n \rightarrow u'$ weakly in $L^{m+1}(\mathbb T^3 \times (0,T))$, then, for each fixed $j$, 
\begin{align}   \label{pre-bre}
\lim_{n\rightarrow \infty} \int_0^t \int_{\mathbb T^3} (|u_j'|^{m-1} u_j')  P_j u_n'  dx d\t
&=  \lim_{n\rightarrow \infty}   \int_0^t \int_{\mathbb T^3} \left[P_j (|u_j'|^{m-1} u_j') \right] u_n'  dx d\t    \notag\\
&= \int_0^t \int_{\mathbb T^3} \left[P_j (|u_j'|^{m-1} u_j') \right] u'  dx d\tau,
\end{align}
for $0\leq t\leq T$. Notice
\begin{align}      \label{bre}
& \int_0^t \int_{\mathbb T^3} \left[P_j (|u_j'|^{m-1} u_j') \right] u'  dx d\tau =  \int_0^t \int_{\mathbb T^3} (|u_j'|^{m-1} u_j') P_j u' dx d\t     \notag\\
&=  \int_0^t \int_{\mathbb T^3} (|u_j'|^{m-1} u_j')  (P_j u'-u') dx d\t + \int_0^t \int_{\mathbb T^3} (|u_j'|^{m-1} u_j')  u' dx d\t.
\end{align}
Since $u'\in L^{m+1} (\mathbb T^3 \times (0,T))$ and $|u_j'|^{m-1} u_j'\rightarrow \psi$ weakly in $L^{\frac{m+1}{m}}(\mathbb T^3 \times (0,T))$, one has
\begin{align}  \label{bre-1}
\lim_{j\rightarrow \infty} \int_0^t \int_{\mathbb T^3} (|u_j'|^{m-1} u_j')  u' dx d\t  = \int_0^t \int_{\mathbb T^3} \psi  u' dx d\t,
\end{align}
for $0\leq t\leq T$. Now we deal with the first integral on the right-hand side of (\ref{bre}). Since $u' \in L^{m+1}(\mathbb T^3 \times (0,T))$, 
$\lim_{j\rightarrow \infty} \|P_j u' - u'\|_{m+1} =0$ for a.e. $t\in [0,T]$, due to (\ref{PSC}). 
Also, $\|P_j u'-u'\|_{m+1} \leq \|P_j u'\|_{m+1}  + \|u'\|_{m+1} \leq c_m\|u'\|_{m+1}$ by (\ref{PSC1}). 
Hence, applying the Lebesgue's dominated convergence theorem, one has  
$\lim_{j\rightarrow \infty} \int_0^T \|P_j u'-u'\|_{m+1}^{m+1} dt =0$. Then, employing the H\"older's inequality,
\begin{align}   \label{bre-2}
\int_0^T \int_{\mathbb T^3} |u_j'|^m  |P_j u'-u'|  dx dt  &\leq   \|u_j'\|_{L^{m+1}(\mathbb T^3 \times (0,T))}       \left(\int_0^T \|P_j u'-u'\|_{m+1}^{m+1} dt \right)^{\frac{1}{m+1}}   \notag\\
&\leq    C   \left( \int_0^T \|P_j u'-u'\|_{m+1}^{m+1} dt \right)^{\frac{1}{m+1}}    \longrightarrow 0,      \text{\;\;as\;\;} j\rightarrow \infty,
\end{align}
where we have used the fact that $u_j'$ is uniformly bounded in $L^{m+1}(\mathbb T^3 \times (0,T))$ by (\ref{ub-4}).

Combining (\ref{bre}), (\ref{bre-1}) and (\ref{bre-2}) gives
\begin{align} \label{bre-3}
\lim_{j\rightarrow \infty}\int_0^t \int_{\mathbb T^3} \left[P_j (|u_j'|^{m-1} u_j') \right] u'  dx d\tau =\int_0^t \int_{\mathbb T^3} \psi  u' dx d\t.
\end{align}
Then, by (\ref{pre-bre}) and (\ref{bre-3}), one has
\begin{align} \label{bre-4}
\lim_{j\rightarrow \infty} \lim_{n\rightarrow \infty}  \int_0^t  \int_{\mathbb T^3} (|u_j'|^{m-1} u_j')  P_j u_n'  dx d\t = \int_0^t \int_{\mathbb T^3} \psi u' dx d\t, \;\; \text{for all} \;\;  t\in [0,T].
\end{align}
Finally, by (\ref{diff-damp}), (\ref{ppre-bre}) and (\ref{bre-4}), we conclude
\begin{align}  \label{conv-junk}
\lim_{j\rightarrow \infty} \lim_{n\rightarrow \infty} \int_0^t \int_{\mathbb T^3} \left[(P_n-P_j)      (|u_j'|^{m-1} u_j' )  \right] w' dx d\t =0,   \;\;  \text{for all} \;\;  t\in [0,T].
\end{align}

In the same manner, we can show 
\begin{align}    \label{conv-junk1}
\lim_{j\rightarrow \infty} \lim_{n\rightarrow \infty} \int_0^t \int_{\mathbb T^3} [(P_n-P_j) (|u_j|^{p-1} u_j) ] w' dx d\tau=0,    \;\;  \text{for all} \;\;  t\in [0,T].
\end{align}

\subsubsection{Estimate for the integral $\int_0^t \int_{\mathbb T^3} \left(|u_n|^{p-1} u_n  - |u_j|^{p-1}  u_j \right)   w' dx d\t$}      \label{AB}

\noindent

\textbf{Case A}:  $1\leq p< \frac{5}{6}(m+1)$.

We estimate
\begin{align}    \label{case1}
&\int_0^T \int_{\mathbb T^3}  \left||u_n|^{p-1} u_n  -  |u_j|^{p-1} u_j \right|  |w'| dx dt   \notag\\
&\leq   C\int_0^T    \int_{\mathbb T^3}  |u_n-u_j| (|u_n|^{p-1}   +  |u_j|^{p-1}) |w'| dx dt    \notag\\
&\leq    C   \|u_n - u_j\|_{  L^{   \frac{m+1}{m+1-p}     }(\mathbb T^3 \times (0,T))}     \left(  \|u_n\|_{L^{m+1}(\mathbb T^3 \times (0,T))}^{p-1}      +      \|u_j\|_{L^{m+1}(\mathbb T^3 \times (0,T))} ^{p-1}   \right)              \|w'\|_{L^{m+1}(\mathbb T^3 \times (0,T))},
\end{align}
where we use the H\"older's inequality in the last step.

Notice that $p< \frac{5}{6}(m+1)$ implies $\frac{m+1}{m+1-p}<6$. 
So there exists $\epsilon_0 \in (0,1)$ such that  $H^{1-\epsilon_0}(\mathbb T^3) \hookrightarrow L^{\frac{m+1}{m+1-p}}(\mathbb T^3)$.
Then, by (\ref{strong}), there is a subsequence $u_n \rightarrow u$ in $C([0,T]; H^{1-\epsilon_0}(\mathbb T^3))$, which implies that $u_n \rightarrow u$ in $L^{   \frac{m+1}{m+1-p}     }(\mathbb T^3 \times (0,T))$.  It follows that
\begin{align}   \label{case1-1}
\lim_{n, j\rightarrow \infty}    \|u_n - u_j\|_{L^{   \frac{m+1}{m+1-p}     }(\mathbb T^3 \times (0,T))}=0.
\end{align}

Since $w'=u_n' - u_j'$ and $u_n'$ is uniformly bounded in $L^{m+1}(\mathbb T^3 \times (0,T))$,
one has $w'$ is uniformly bounded in $L^{m+1}(\mathbb T^3 \times (0,T))$.
Also, recall $u_n$ is uniformly bounded in $L^{\infty}(0,T; L^{m+1}(\mathbb T^3))$ by (\ref{ub-2}).
Then, by (\ref{case1}) and (\ref{case1-1}), we have
\begin{align}    \label{fCauchy}
\lim_{n, j\rightarrow \infty}  \int_0^T \int_{\mathbb T^3}  \left||u_n|^{p-1} u_n  -  |u_j|^{p-1} u_j \right|  |w'| dx dt  =0.
\end{align}

Applying (\ref{conv-junk}), (\ref{conv-junk1}) and (\ref{fCauchy}) to  inequality (\ref{diff-ener1}),
and since $\lim_{n, j \rightarrow \infty}  (\|\nabla w(0)\|_2^2 + \|w'(0)\|_2^2) = \lim_{n,j\rightarrow \infty}   (\|\nabla (P_n -   P_j) u_0\|_2^2 + \|(P_n-P_j)  u_1\|_2^2) =   0$, we obtain 
\begin{align}     \label{gCauchy}
\lim_{j \rightarrow \infty}   \varlimsup_{n\rightarrow \infty}  \int_0^T \int_{\mathbb T^3} \left( |u_n'|^{m-1} u_n'  -   |u_j'|^{m-1} u_j  \right)w' dx d\t  =0. 
\end{align}

\begin{remark}
Recall Case 1 of the assumption for the ``existence" result (Theorem \ref{thm1}) is that $1\leq p\leq m \leq 5$. Notice
\begin{align}
\{(m,p) \not = (5,5):  \, 1\leq p\leq m \leq 5    \}  \subset  \{(m,p):   1\leq p< \frac{5}{6}(m+1)  \}.
\end{align}
Nonetheless, the situation $m=p=5$ does not satisfy $1\leq p< \frac{5}{6}(m+1)$, thus we have to discuss it separately.
\end{remark}

\textbf{Case B}: $m=p=5$.  

In this case,  $$\int_0^t \int_{\mathbb T^3} \left(|u_n|^{p-1} u_n  - |u_j|^{p-1}  u_j \right)   w' dx d\tau =   \int_0^t \int_{\mathbb T^3} 
\left( u_n^5  -  u_j^5 \right)   w' dx d\tau .$$
We estimate the above integral by using integration by parts in time. Such ``integration by parts" technique originates from \cite{BL1}.
Indeed, since $w=u_n-u_j$,
\begin{align}  \label{mp5}
&\left|\int_0^t \int_{\mathbb T^3}  \left(u_n^5- u_j^5 \right) w' dx d\tau  \right| \notag\\
&= \frac{1}{2} \left|\int_0^t  \int_{\mathbb T^3}  (u_n^4 +  u_n^3 u_j + u_n^2 u_j^2 + u_n u_j^3 + u_j^4)   (w^2)' dx d\tau \right|  \notag\\
&\leq  C\int_{\mathbb T^3}   \left(|u_n(t)|^4  +  |u_j(t)|^4\right) |w(t)|^2 dx +    C\int_{\mathbb T^3}   \left(|u_n(0)|^4  +  |u_j(0)|^4\right) |w(0)|^2 dx  \notag\\
& \hspace{0.3 in}  +    C \int_0^t \int_{\mathbb T^3}  \left(|u_n|^3+ |u_j|^3\right)    \left(|u_n'|+ |u_j'|\right)   w^2   dx d\tau.
\end{align}
We shall estimate each term on the right-hand side of (\ref{mp5}). 

Since $m=5$, then owing to (\ref{ub-1})-(\ref{ub-4}), there exists a uniform bound $K$ such that
\begin{align}    \label{mp5-0}
&\|u_n(t)\|_{H^1}^2 +   \|u_n'(t)\|_2^2 +  \int_0^t  \|u_n'(\tau)\|_6^6 d\tau   \leq K, \;\; \text{for all} \,\, t\in [0,T], \;\text{for all} \,\, n\in \mathbb N.
\end{align}

By H\"older's inequality, 
\begin{align}  \label{mp5-1}
&\int_0^t \int_{\mathbb T^3}  \left(|u_n|^3+ |u_j|^3\right)    \left(|u_n'|+ |u_j'|\right)   w^2   dx d\tau  \leq   C  \int_0^t   \left(\|u_n\|_6^3   + \|u_j\|_6^3 \right)  \left(\|u_n'\|_6 +   \|u_j'\|_6  \right)    \|w\|_6^2      d\tau   \notag\\
&\leq   C(K) \int_0^t   \left(\|u_n'\|_6 +   \|u_j'\|_6  \right)    \|w\|_{H^1}^2      d\tau,   \;\;   \text{for all} \,\, t\in [0,T],
\end{align}
where we use $H^1 \hookrightarrow L^6$ as well as (\ref{mp5-0}) to obtain the last inequality.

Also, using H\"older's inequality, we have
\begin{align}    \label{mp5-2}
&\int_{\mathbb T^3}   \left(|u_n(0)|^4  +  |u_j(0)|^4\right) |w(0)|^2 dx    \notag\\
&\leq    \left(\|u_n(0)\|_6^4   +    \|u_j(0)\|_6^4\right)    \|w(0)\|_6^2   \leq    C(K)    \|w(0)\|_{H^1}^2,
\end{align}
where the last inequality is due to (\ref{mp5-0}).

Next, we estimate $\int_{\mathbb T^3}   \left(|u_n(t)|^4  +  |u_j(t)|^4\right) |w(t)|^2 dx$. Indeed,
\begin{align}   \label{mp5-3}
\int_{\mathbb T^3}    |u_n(t)|^4   |w(t)|^2 dx     \leq    & C\int_{\mathbb T^3}    |u_n(t) - u_n(0)|^4   |w(t)|^2 dx    
+   C \int_{\mathbb T^3}    |u_n(0)  - u_0|^4   |w(t)|^2 dx \notag\\
&  +  C \int_{\mathbb T^3}    |u_0|^4   |w(t)|^2 dx.
\end{align}
We estimate each term on the right-hand side of (\ref{mp5-3}). First,
\begin{align}   \label{mp5-4}
&\int_{\mathbb T^3}    |u_n(t) - u_n(0)|^4   |w(t)|^2 dx  =
\int_{\mathbb T^3}    \left|  \int_0^t   u_n'(\tau)  d\tau \right|^4       |w(t)|^2  dx     \notag\\
&\leq   \left(\int_{\mathbb T^3}     \left|  \int_0^t   u_n'(\tau)  d\tau \right|^6  dx \right)^{\frac{2}{3}}      \|w(t)\|_6^2     \notag\\
&\leq C t^{\frac{10}{3}} \left( \int_0^t   \|u_n'(\tau)\|_6^6 d\tau  \right)^{\frac{2}{3}}   \|w(t)\|_{H^1}^2 
\leq C(K) t^{\frac{10}{3}}  \|w(t)\|_{H^1}^2,   \;\;  \text{for all} \,\,  t\in [0,T],
\end{align}
where the last inequality is due to (\ref{mp5-0}). Also, since $u_n(0)  = P_n u_0 \rightarrow u_0$ in $H^1$, then
\begin{align}     \label{mp5-5}
&\int_{\mathbb T^3}    |u_n(0)  - u_0|^4   |w(t)|^2 dx  
\leq \|u_n(0) - u_0\|_6^4    \|w(t)\|_6^2    \notag\\
&\leq C  \|u_n(0) - u_0\|_{H^1}^4   \|w(t)\|_{H^1}^2 \leq \epsilon \|w(t)\|_{H^1}^2, \;\;  \text{for} \,\,  n \geq N_{\epsilon}.
\end{align}
Notice 
\begin{align}   \label{mp5-66}
\|w(t)\|_2^2 = \int_{\mathbb T^3} |w(t)|^2 dx 
&\leq   \int_{\mathbb T^3}  \left|\int_0^t w'(\tau) d\tau\right|^2 dx + \|w(0)\|_2^2   \notag\\
&\leq     t  \int_0^t \|w'(\tau)\|_2^2 d\tau    +   \|w(0)\|_2^2.
\end{align}
We let $\varphi$ be a periodic smooth function such that $\|u_0 - \varphi\|_{H^1}^4 \leq \epsilon$. 
Since $\varphi$ is smooth, there exists $C_{\epsilon}>0$ such that $|\varphi(x)| \leq C_{\epsilon}$ for all $x\in \mathbb T^3$. As a result,
\begin{align}     \label{mp5-6}
&\int_{\mathbb T^3}    |u_0|^4   |w(t)|^2 dx
\leq    C\int_{\mathbb T^3}    |u_0 - \varphi|^4   |w(t)|^2 dx   +     C\int_{\mathbb T^3}    |\varphi|^4   |w(t)|^2 dx    \notag\\
&\leq  C\|u_0-\varphi\|_6^4   \|w(t)\|_6^2  +   C_{\epsilon} \|w(t)\|_2^2 \notag\\
&\leq     C\epsilon \|w(t)\|_{H^1}^2  +    C_{\epsilon}   t  \int_0^t \|w'(\tau)\|_2^2 d\tau    +   C_{\epsilon}   \|w(0)\|_2^2,
\;\;\text{for all}\,\, t\in [0,T],
\end{align}
owing to (\ref{mp5-66}).

Applying estimates (\ref{mp5-4}), (\ref{mp5-5}) and (\ref{mp5-6}) to the right-hand side of (\ref{mp5-3}) yields
\begin{align}    \label{mp5-7}
\int_{\mathbb T^3}    |u_n(t)|^4   |w(t)|^2 dx   \leq    C(K) ( t^{\frac{10}{3}} + \epsilon) \|w(t)\|_{H^1}^2 
+   C_{\epsilon}  t \int_0^t \|w'(\tau)\|_2^2 d\tau    +   C_{\epsilon}   \|w(0)\|_2^2,
\end{align}
for all $t\in [0,T]$, and for all $n\geq N_{\epsilon}$. Thus, $\int_{\mathbb T^3}   \left(|u_n(t)|^4  +  |u_j(t)|^4\right) |w(t)|^2 dx$ is also bounded by the right-hand side of (\ref{mp5-7}) if $n, j\geq N_{\epsilon}$.

By virtue of (\ref{mp5}), (\ref{mp5-1}), (\ref{mp5-2}) and (\ref{mp5-7}), we obtain
\begin{align}    \label{mp5-8}
&\left| \int_0^t \int_{\mathbb T^3}  \left(u_n^5- u_j^5 \right) w' dx d\tau  \right| \notag\\
&\leq     C(K) ( t^{\frac{10}{3}} + \epsilon)\|w(t)\|_{H^1}^2  +   C_{\epsilon}   t  \int_0^t \|w'(\tau)\|_2^2 d\tau    \notag\\
& \hspace{0.2 in } + C(K) \int_0^t   \left(\|u_n'\|_6 +   \|u_j'\|_6  \right)    \|w\|_{H^1}^2      d\tau 
+   C(K,\epsilon)   \|w(0)\|_{H^1}^2,  
\end{align}
for all $t\in [0,T]$ and $n, j\geq N_{\epsilon}$.

Applying estimates (\ref{mp5-66}) and (\ref{mp5-8}) to inequality (\ref{diff-ener1}) with $m=p=5$, 
and using $\|w\|^2_{H^1} =  \|\nabla w\|_2^2 +  \|w\|_2^2$, we have 
\begin{align}    \label{mp5-9}
 0&\leq \frac{1}{2}\left(   \|w(t)\|^2_{H^1} + \|w'(t)\|^2_2 \right)  +        \int_0^t \int_{\mathbb T^3} \left[ (u_n')^5- (u_j')^5 \right]  w' dx d\t  \notag\\
& \leq      C(K,\epsilon) \|w(0)\|^2_{H^1}   +   \frac{1}{2} \|w'(0)\|^2_2   
 +  C(K) ( t^{\frac{10}{3}}  + \epsilon) \|w(t)\|_{H^1}^2     \notag\\
& \hspace{0.2 in }   +   C_{\epsilon}   t \int_0^t \|w'(\tau)\|_2^2 d\tau    + C(K)  \int_0^t   \left(\|u_n'\|_6 +   \|u_j'\|_6  \right)    \|w\|_{H^1}^2      d\tau
    \notag\\
&  \hspace{0.2 in} + \left|\int_0^t \int_{\mathbb T^3}  \left[(P_n -P_j)    u_j^5   \right]  w' dx d\tau \right|
+  \left|\int_0^t \int_{\mathbb T^3} \left[(P_n-P_j)   (u_j')^5 \right]   w' dx d\tau\right|,
\end{align}
for all $t\in [0,T]$ and $n, j\geq N_{\epsilon}$.

We remark that our strategy is to first prove the local existence of weak solutions on $[0,T]$, and extend the local solution to a global solution later.
Thus, we can choose $\epsilon$ and $T$ sufficiently small 
such that $C(K) (T^{\frac{10}{3}} + \epsilon)\leq \frac{1}{4}$, 
then (\ref{mp5-9}) shows
\begin{align}    \label{mp5-10}
0 &\leq \frac{1}{4}\left(   \|w(t)\|^2_{H^1} + \|w'(t)\|^2_2 \right)  +        \int_0^t \int_{\mathbb T^3} \left[ (u_n')^5- (u_j')^5 \right]   w' dx d\t  \notag\\
& \leq      C(K,\epsilon)  \|w(0)\|^2_{H^1}   +   \frac{1}{2} \|w'(0)\|^2_2   + C_{\epsilon}   t \int_0^t \|w'(\tau)\|_2^2 d\tau + C(K)  \int_0^t   \left(\|u_n'\|_6 +   \|u_j'\|_6  \right)    \|w\|_{H^1}^2      d\tau
    \notag\\
&  \hspace{0.2 in} + \left|\int_0^t \int_{\mathbb T^3}  \left[(P_n -P_j)    (u_j^5)   \right]  w' dx d\tau \right|
+  \left|\int_0^t \int_{\mathbb T^3} \left[(P_n-P_j)   (u_j')^5 \right]   w' dx d\tau\right|,
\end{align}      
for all $t\in [0,T]$ and $n, j\geq N_{\epsilon}$.

Because of (\ref{mp5-0}), we can apply the Gr\"onwall's inequality to (\ref{mp5-10}), it follows that
\begin{align}   \label{mp5-11}
0 &\leq \frac{1}{4}\left(   \|w(t)\|^2_{H^1} + \|w'(t)\|^2_2 \right)  +        \int_0^t \int_{\mathbb T^3} \left[ (u_n')^5- (u_j')^5 \right]   w' dx d\t \notag\\
& \leq C(K,T,\epsilon) \Big(\|w(0)\|_{H^1}^2 +  \|w'(0)\|_2^2 +   \left|\int_0^t \int_{\mathbb T^3} [(P_n -P_j) u_j^5 ]w' dx d\t  \right|  \notag\\
&\hspace{1 in}+   \left|\int_0^t \int_{\mathbb T^3} [(P_n-P_j) (u_j')^5 ] w' dx d\tau  \right|  \Big),
\end{align}
for $t\in [0,T]$ and $n, j\geq N_{\epsilon}$.

Since $w=u_n - u_j$, one has 
$$\lim_{n,j\rightarrow \infty} (\|w(0)\|_{H^1}^2 +  \|w'(0)\|_2^2) =  \lim_{n,j\rightarrow \infty} (\|P_n u_0 - P_j u_0\|_{H^1}^2 +  \|P_n u_1 - P_j u_1\|_2^2)=0.$$
Then, by using (\ref{conv-junk}) and (\ref{conv-junk1}), we derive from (\ref{mp5-11}) that
\begin{align*}
\lim_{j\rightarrow \infty} \varlimsup_{n\rightarrow \infty}     \int_0^t \int_{\mathbb T^3} \left[(u_n')^5 - (u_j')^5 \right]w ' dx d\tau  =0,  \;\; \text{for all}\,\, t\in [0,T],
\end{align*}
under the scenario that $m=p=5$.

In sum, for both Case A (i.e. $1\leq p< \frac{5}{6}(m+1)$) and Case B (i.e. $m=p=5$), we have
\begin{align}  \label{fina}
\lim_{j \rightarrow \infty}   \varlimsup_{n\rightarrow \infty}  \int_0^T \int_{\mathbb T^3} \left( |u_n'|^{m-1} u_n'  -   |u_j'|^{m-1} u_j  \right)w' dx d\t  =0.\end{align}
Although in Case B we restrict $T$ to be small, this restriction does not affect our intention to prove the local existence of weak solutions on $[0,T]$. Local weak solutions will eventually be extended to global ones in subsection \ref{sec-exten}.

We remark that the ``union" of Case A and Case B in the above proof is same as the the ``union" of Case 1 and Case 2 in the statement of Theorem \ref{thm1}.

\vspace{0.1 in}

\subsubsection{Completion of the proof for $|u'_n|^{m-1}u_n'  \rightarrow |u'|^{m-1} u'$ in $L^{\frac{m+1}{m}}(\mathbb T^3  \times (0,T))$}
We recall from (\ref{con-4}) and (\ref{g-conv}) that $u'_n \rightarrow u'$ weakly in $L^{m+1}(\mathbb T^3 \times (0,T))$ 
and $|u_n'|^{m-1}u_n' \rightarrow \psi$ weakly in $L^{\frac{m+1}{m}}(\mathbb T^3 \times (0,T))$. Then, by (\ref{fina}), we have 
\begin{align*} 
&\varlimsup_{j\rightarrow \infty} \varlimsup_{n\rightarrow \infty}\int_0^T \int_{\mathbb T^3} (|u_n'|^{m-1}u_n'  -  |u_j'|^{m-1} u_j'  )(u_n'  - u_j')   dx dt \notag\\
&=  \varlimsup_{n\rightarrow \infty}    \int_0^T \int_{\mathbb T^3}   |u_n'|^{m+1} dx dt   +  
\varlimsup_{j\rightarrow \infty}    \int_0^T \int_{\mathbb T^3}    |u_j'|^{m+1} dx dt   \notag\\
& \hspace{0.2 in}  - \lim_{j\rightarrow \infty} \lim_{n\rightarrow \infty} \int_0^T \int_{\mathbb T^3}  (|u_n'|^{m-1}u_n')u_j'  dx dt
 -   \lim_{j\rightarrow \infty} \lim_{n\rightarrow \infty} \int_0^T \int_{\mathbb T^3}  (|u_j'|^{m-1} u_j) u_n'  dx dt \notag\\
&=  2\varlimsup_{n\rightarrow \infty}    \int_0^T \int_{\mathbb T^3}   |u_n'|^{m+1} dx dt   - 2   \int_0^T \int_{\mathbb T^3}  \psi  u'  dx dt =0.
\end{align*}
It follows that
\begin{align}    \label{lsup}
\varlimsup_{n\rightarrow \infty}    \int_0^T \int_{\mathbb T^3}    |u_n'|^{m+1}  dx dt
=  \int_0^T \int_{\mathbb T^3}  \psi  u'  dx dt.
\end{align}
Since $|s|^{m-1}s$ is a monotone increasing function on $\mathbb R$, then for any $v \in L^{m+1}(\mathbb T^3 \times (0,T))$, we have 
\begin{align}   \label{positivity}
\int_0^T  \int_{\mathbb T^3} (|u_n'|^{m-1}u_n'  - |v|^{m-1}v ) (u_n' - v) dx dt  \geq 0.
\end{align}
Owing to (\ref{lsup}) and (\ref{positivity}), as well as the fact that $u'_n \rightarrow u'$ weakly in $L^{m+1}(\mathbb T^3 \times (0,T))$ 
and $|u_n'|^{m-1}u_n' \rightarrow \psi$ weakly in $L^{\frac{m+1}{m}}(\mathbb T^3 \times (0,T))$, we obtain
\begin{align}   \label{mm}
&\varlimsup_{n\rightarrow \infty}   \int_0^T  \int_{\mathbb T^3} (|u_n'|^{m-1}u_n'  - |v|^{m-1}v ) (u_n' - v) dx dt \notag\\
&=  \varlimsup_{n\rightarrow \infty}   \int_0^T  \int_{\mathbb T^3}   |u_n'|^{m+1}   dx dt
-     \lim_{n\rightarrow \infty}   \int_0^T    \int_{\mathbb T^3}    (|u_n'|^{m-1}u_n')  v  dx dt   \notag\\
& \hspace{0.2 in}  -    \lim_{n\rightarrow \infty}   \int_0^T    \int_{\mathbb T^3}    (|v|^{m-1}v) u_n'  dx dt 
 +    \int_0^T    \int_{\mathbb T^3}   |v|^{m+1}  dx dt    \notag\\
&=     \int_0^T    \int_{\mathbb T^3}    \psi  u' d x dt -       \int_0^T    \int_{\mathbb T^3}    \psi v  d x dt  
 -   \int_0^T   \int_{\mathbb T^3} (|v|^{m-1}v)  u' dx dt +     \int_0^T    \int_{\mathbb T^3}   |v|^{m+1} dx dt \notag\\
& =    \int_0^T    \int_{\mathbb T^3}    (\psi   - |v|^{m-1}v)(u' -v) dx dt \geq 0,  \;\; \text{for any} \,\, v \in L^{m+1}(\mathbb T^3 \times (0,T)).
\end{align}
We claim that $v\mapsto |v|^{m-1}v$ is a maximal monotone operator mapping from $L^{m+1}(\mathbb T^3 \times (0,T))  \rightarrow L^{\frac{m+1}{m}}(\mathbb T^3 \times (0,T))$. It suffices to show that this operator is both monotone and hemicontinuous (see, e.g., \cite{Ba}). The monotonicity is obvious. It remains to verify the hemicontinuity. Recall that an operator $A$ mapping from a Banach space to its dual is said to be hemicontinuous if $A(v+\lambda y)$ converges weakly to $A(v)$ as $\lambda \rightarrow 0$ for all $v$, $y$ in this Banach space. Hence, to check the operator $v\mapsto |v|^{m-1}v$ is hemicontinuous from $L^{m+1}(\mathbb T^3 \times (0,T))  \rightarrow L^{\frac{m+1}{m}}(\mathbb T^3 \times (0,T))$, it is enough to verify that, for all $v$, $y$, $\eta \in   L^{m+1}(\mathbb T^3 \times (0,T)) $,
\begin{align}    \label{hemi}
\lim_{\lambda\rightarrow \infty}  \int_0^T \int_{\mathbb T^3}   \left[|v +  \lambda y |^{m-1} (v + \lambda y) \right] \eta dx   dt  =  \int_0^T  \int_{\mathbb T^3} (|v|^{m-1} v) \eta  dx dt.
\end{align}
As a matter of fact, (\ref{hemi}) follows from Lebesgue's dominated convergence theorem. 

Since we have shown the maximal monotonicity of the operator $v\mapsto |v|^{m-1}v$ from $L^{m+1}(\mathbb T^3 \times (0,T))  \rightarrow L^{\frac{m+1}{m}}(\mathbb T^3 \times (0,T))$, it can be concluded from (\ref{mm}) that $\psi = |u'|^{m-1} u'$. Namely,
\begin{align}   \label{conv-gn}
|u_n'|^{m-1}u_n' \rightarrow |u'|^{m-1} u'  \;\; \text{weakly in} \,\, L^{\frac{m+1}{m}}(\mathbb T^3 \times (0,T)).
\end{align}

\vspace{0.1 in}

\subsection{Passage to the limit for the Galerkin system}
In this section, we let $n\rightarrow \infty$ in the Galerkin approximation system, and aim to show that $(u,u')$ is a weak solution for model (\ref{sys1})-(\ref{sys2}) on $[0,T]$. Let $\phi$ be a trigonometric polynomial with smooth coefficients, i.e., 
$\phi=\sum_{\substack{k=(k_1,k_2,k_3) \in \mathbb Z^3\\ |k_1|, |k_2|, |k_3| \leq N}} \hat{\phi}(k,t) e^{ik\cdot x}$
where $\hat{\phi}(k,t)$ is smooth in $t$. We multiply the Galerkin system (\ref{Gal}) with $\phi$ and integrate over $\mathbb T^3 \times (0,t)$. Assume $n$ is larger than the degree $N$ of $\phi$, then 
$\int_0^t \int_{\mathbb T^3}  P_n(|u_n'|^{m-1} u_n') \phi    dx d\tau =   \int_0^t \int_{\mathbb T^3} (|u_n'|^{m-1} u_n') \phi    dx d\tau$,
and $\int_0^t \int_{\mathbb T^3}  P_n(|u_n|^{p-1} u_n) \phi    dx d\tau =   \int_0^t \int_{\mathbb T^3} (|u_n|^{p-1} u_n) \phi    dx d\tau$. 
It follows that
\begin{align}    \label{bpass}
&  \int_0^t \int_{\mathbb T^3} u_n''  \phi  dx d\tau   +  \int_0^t  \int_{\mathbb T^3} \nabla u_n  \cdot \nabla  \phi dx d\tau +  
  \int_0^t \int_{\mathbb T^3}   (|u_n'|^{m-1} u_n')   \phi dx d\tau    \notag\\
 & = \int_0^t   \int_{\mathbb T^3}  (|u_n|^{p-1} u_n) \phi dx d\tau,   \;\;   \text{for all}    \;\;  t\in [0,T].
\end{align}
Then, since $|u_n'|^{m-1}u_n' \rightarrow |u'|^{m-1} u'$ weakly in $L^{\frac{m+1}{m}}(\mathbb T^3 \times (0,T))$, $|u_n|^{p-1}u_n    \rightarrow |u|^{p-1}u $ weakly in $L^{\frac{m+1}{m}}(\mathbb T^3 \times (0,T))$, 
$u_n \rightarrow u$ weakly$^*$ in $L^{\infty}(0,T;  H^1(\mathbb T^3))$ 
and $u_n'' \rightarrow u''$ weakly$^*$ in $L^{\frac{m+1}{m}}(0,T;  X')$ where $X=H^1(\mathbb T^3) \cap L^{m+1}(\mathbb T^3)$, we can pass to the limit as $n\rightarrow \infty$ in (\ref{bpass}) to obtain 
\begin{align}    \label{var}
&  \int_0^t  \langle u'',  \phi \rangle d\tau   +  \int_0^t  \int_{\mathbb T^3} \nabla u \cdot \nabla  \phi dx d\tau +  
  \int_0^t \int_{\mathbb T^3}   (|u'|^{m-1} u')   \phi dx d\tau    \notag\\
 & = \int_0^t   \int_{\mathbb T^3}  (|u|^{p-1} u) \phi dx d\tau,   \;\;   \text{for all}    \;\;  t\in [0,T].
\end{align}
After integration by parts in time, we obtain
\begin{align}  \label{var2}   
&\int_{\mathbb T^3} u'(t) \phi(t) dx - \int_{\mathbb T^3} u'(0) \phi(0) dx - \int_0^t \int_{\mathbb T^3} u'  \phi'  dx d\tau   \notag\\
&+  \int_0^t  \int_{\mathbb T^3} \nabla u   \cdot \nabla  \phi dx d\tau +  
  \int_0^t \int_{\mathbb T^3}   (|u'|^{m-1} u' ) \phi dx d\tau = \int_0^t   \int_{\mathbb T^3} (|u|^{p-1} u) \phi dx d\tau,
\end{align}
for all $t\in [0,T]$, and for any trigonometric polynomial $\phi$ with smooth coefficients.

Recall $u \in L^{\infty}(0,T;H^1(\mathbb T^3))$, $u' \in L^{\infty}(0,T;L^2(\mathbb T^3))$, and
$|u|^{p-1}u$, $|u'|^{m-1}u' \in L^{\frac{m+1}{m}}(\mathbb T^3 \times (0,T))$. Thus, by density, we conclude that (\ref{var2}) holds for all $\phi \in C([0,T]; H^1(\mathbb T^3)) \cap L^{m+1}(\mathbb T^3 \times (0,T))$ with $\phi' \in C([0,T];L^2(\mathbb T^3))$.

We shall verify the initial condition $u(0)=u_0$ and $u'(0)=u_1$. Indeed, let $\tilde \phi$ be a trigonometric polynomial with smooth coefficients such that $\tilde \phi(T)={\tilde \phi}'(T)=0$. By setting $t=T$ and $\phi=\tilde \phi$ in (\ref{bpass}), then after integration by parts in time, we obtain 
\begin{align}   \label{IC-1}
& \int_{\mathbb T^3} [    (P_n u_0) \tilde \phi'(0)   -      (P_n u_1) \tilde \phi(0)     ] dx + \int_0^T \int_{\mathbb T^3} u_n  \tilde \phi'' dx d\tau   +  \int_0^T \int_{\mathbb T^3} \nabla u_n  \cdot \nabla  \tilde \phi dx d\tau  \notag\\  
  & +\int_0^T \int_{\mathbb T^3}   (|u_n'|^{m-1} u_n')  \tilde \phi dx d\tau   
  = \int_0^T   \int_{\mathbb T^3}  (|u_n|^{p-1} u_n) \tilde \phi dx d\tau,
\end{align}
where we have used $u_n(0)=P_n u_0$ and $u_n'(0)=P_n u_1$.

Letting $n\rightarrow \infty$ in (\ref{IC-1}), we have
\begin{align}   \label{IC-2}
& \int_{\mathbb T^3} [    u_0 \tilde \phi'(0)   -     u_1\tilde \phi(0)     ] dx + \int_0^T \int_{\mathbb T^3} u \tilde \phi'' dx d\tau   +  \int_0^T \int_{\mathbb T^3} \nabla u  \cdot \nabla  \tilde \phi dx d\tau  \notag\\  
  & +\int_0^T \int_{\mathbb T^3}   (|u'|^{m-1} u')  \tilde \phi dx d\tau   
  = \int_0^T   \int_{\mathbb T^3}  (|u|^{p-1} u) \tilde \phi dx d\tau.
\end{align}
However, by setting $t=T$ and $\phi=\tilde \phi$ in (\ref{var}), and after integration by parts, we have
\begin{align}   \label{IC-3}
& \int_{\mathbb T^3} [    u(0) \tilde \phi'(0)   -     u'(0) \tilde \phi(0)     ] dx + \int_0^T \int_{\mathbb T^3} u  \tilde \phi'' dx d\tau   +  \int_0^T \int_{\mathbb T^3} \nabla u  \cdot \nabla  \tilde \phi dx d\tau  \notag\\  
  & +\int_0^T \int_{\mathbb T^3}   (|u'|^{m-1} u')   \tilde \phi dx d\tau   
  = \int_0^T   \int_{\mathbb T^3}  (|u|^{p-1} u)\tilde \phi dx d\tau.
\end{align}
Comparing (\ref{IC-2}) and (\ref{IC-3}), we obtain   
\begin{align}   \label{IC-4}
\int_{\mathbb T^3} [    u_0\tilde \phi'(0)   -     u_1 \tilde \phi(0)     ] dx  =   \int_{\mathbb T^3} [    u(0) \tilde \phi'(0)   -     u'(0)\tilde \phi(0)     ] dx.
\end{align}
Since $\tilde \phi(0)$ and $\tilde \phi'(0)$ are arbitrary trigonometric polynomials, we obtain $u(0)=u_0$ and $u'(0)=u_1$. 
This completes the proof for the existence of weak solutions on $[0,T]$.

\vspace{0.1 in}

\subsection{Energy identity}   \label{sec-EI}
Let $(u,u_t)$ be a weak solution for (\ref{sys1})-(\ref{sys2}) on $[0,T]$ in the sense of Definition \ref{weak}. We aim to show that the energy identity (\ref{EI}) holds for all $t\in [0,T]$.

One may multiply equation (\ref{sys1}) with $u_t$, and perform integration by parts to obtain the energy identity (\ref{EI}) \emph{formally}. However, $u_t$ is not smooth enough, this formal procedure is not rigorous. To remedy it, we regularize solutions by the projection operator $P_n$ defined in (\ref{Spartial}). We set $\phi=P_n  u_t$ in the variational identity (\ref{varI}) to get
\begin{align}  \label{pene}   
&\frac{1}{2} \left(\|P_n u_t(t)\|_2^2 + \|\nabla P_n u(t)\|_2^2\right) + \int_0^t \int_{\mathbb T^3}   (|u_t|^{m-1} u_t)  (P_n u_t) dx d\tau \notag\\
&=     \frac{1}{2} \left(\|P_n u_t(0)\|_2^2 + \|\nabla P_n u(0)\|_2^2\right) +   \int_0^t   \int_{\mathbb T^3} (|u|^{p-1} u) (P_n u_t)dx d\tau,
\end{align}
for all $t\in [0,T]$.

Since $u_t \in L^{m+1}(\mathbb T^3 \times (0,T))$, one has $u_t \in L^{m+1}(\mathbb T^3)$ for a.e. $t\in [0,T]$. Then by (\ref{PSC}), $P_n u_t\rightarrow u_t$ in $L^{m+1}(\mathbb T^3)$ as $n\rightarrow \infty$ for a.e. $t\in [0,T]$. Also, due to (\ref{PSC1}), we know $\|P_n u_t\|_{m+1}\leq c_m \|u_t\|_{m+1}$. 
Therefore, by the Lebesgue's dominated convergence theorem, we obtain
\begin{align}   \label{pene0}
\lim_{n\rightarrow \infty} \int_0^t \int_{\mathbb T^3} |P_n u_t- u_t|^{m+1} dx d\tau  =0.
\end{align}

Thus, by H\"older's inequality, 
\begin{align*}
&\left|\int_0^t \int_{\mathbb T^3}  (|u_t|^{m-1} u_t)  (P_n u_t-u_t) dx d\tau \right|  \notag\\
&\leq  \left(\int_0^t \int_{\mathbb T^3} |u_t|^{m+1} dx d\tau\right)^{\frac{m}{m+1}}  \left( \int_0^t \int_{\mathbb T^3} |P_n u_t- u_t|^{m+1} dx d\tau\right)^{\frac{1}{m+1}}
 \rightarrow 0,  \;\;  \text{as}\,\, n\rightarrow \infty,
\end{align*}
where the convergence to zero is due to (\ref{pene0}). It follows that
\begin{align}   \label{pene1}
\lim_{n\rightarrow \infty}  \int_0^t \int_{\mathbb T^3} (|u_t|^{m-1} u_t) (P_n u_t) dx d\tau = \int_0^t \int_{\mathbb T^3} |u_t|^{m+1} dx d\tau,  \;\;  \text{for all}  \;\; t\in [0,T].
\end{align}

Analogously, we can derive
\begin{align}  \label{pene4}
\lim_{n\rightarrow \infty} \int_0^t \int_{\mathbb T^3}   (|u|^{p-1} u)  (P_n u_t) dx d\tau &= \int_0^t \int_{\mathbb T^3}   |u|^{p-1}u u_t dx d\tau,
\;\;\text{for all}\,\, t\in [0,T].
\end{align}
Thanks to (\ref{pene1}) and (\ref{pene4}), we let $n\rightarrow \infty$ in (\ref{pene}) to obtain the energy identity
\begin{align}  \label{pene5}   
&\frac{1}{2} \left(\|u_t(t)\|_2^2 + \|\nabla u(t)\|_2^2\right) + \int_0^t \int_{\mathbb T^3}   |u_t|^{m+1}    dx d\tau \notag\\
&=     \frac{1}{2} \left(\|u_t(0)\|_2^2 + \|\nabla u(0)\|_2^2\right) +   \int_0^t   \int_{\mathbb T^3} |u|^{p-1} u u_t dx d\tau,  
\;\;\text{for all}\,\, t\in [0,T].
\end{align}
Since $ \int_0^t   \int_{\mathbb T^3} |u|^{p-1} u u_t dx d\tau = \frac{1}{p+1}\left(\|u(t)\|_{p+1}^{p+1}-  \|u(0)\|_{p+1}^{p+1}\right)$,
the energy identity (\ref{pene5}) can be written in the form of (\ref{EI}) stated in Theorem \ref{thm1}.

\vspace{0.1 in}

\subsection{Extension to global solutions}     \label{sec-exten}
Let $(u,u_t)$ be a weak solution for (\ref{sys1})-(\ref{sys2}) on $[0,T_{max})$ where $T_{max}\in (0,\infty]$ is the maximum life span of this solution. 
Thus, energy identity (\ref{pene5}) holds for all $t\in [0,T_{max})$. From (\ref{pene5}) and performing the same calculation as in (\ref{ener-idn0})-(\ref{exp-grow}), we can derive 
\begin{align}     \label{pext}
\mathscr E(t)  +   \frac{1}{2} \int_0^t  \|u_t(\tau)\|^{m+1}_{m+1} d\t \leq  \left(\mathscr E(0) + t\right)      e^{Ct}, \;\;   \text{for all}    \,\, t\in [0,T_{max}),
\end{align}
where $\mathscr E(t):=    \frac{1}{2}\left(\|\nabla u(t)\|^2_2+\|u_t(t)\|^2_2 \right)+\frac{1}{m+1}\norm{u(t)}_{m+1}^{m+1}$. 
By virtue of (\ref{pext}), we conclude that the weak solution does not blow up at $T_{max}$. Therefore, $T_{max}=+\infty$. Also, (\ref{pext}) tells us that the weak solution grows at most exponentially in time as $t\rightarrow \infty$. This completes the proof for Theorem \ref{thm1}.

\vspace{0.1 in}

\section{Uniqueness of weak solutions and continuous dependence on initial data}
This section is devoted to proving Theorem \ref{thm2}, namely, the uniqueness of weak solutions as well as the continuous dependence on initial data.
We present the proof for the continuous dependence on initial data. The uniqueness of weak solutions follows immediately.

Suppose $(u_0, u_1) \in \left(H^1(\mathbb T^3) \cap L^{m+1}(\mathbb T^3)\right) \times L^2(\mathbb T^3)$.  Let $(u_0^n,u_1^n)$ be a sequence of periodic functions 
in $\left(H^1(\mathbb T^3) \cap L^{m+1}(\mathbb T^3)\right) \times L^2(\mathbb T^3)$ such that $\lim_{n \rightarrow \infty}\|u_0^n - u_0\|_{H^1}=0$, $\lim_{n \rightarrow \infty}\|u_0^n - u_0\|_{m+1}=0$, and $\lim_{n \rightarrow \infty}\|u_1^n - u_1\|_2=0$.
In the previous section, we have proved the existence of global weak solutions for system (\ref{sys1})-(\ref{sys2}). Thus, for any $T>0$, there exists a weak solution $(u,u')$ for system (\ref{sys1})-(\ref{sys2}) on $[0,T]$ with initial data $(u_0,u_1)$. 
Also, for each $n\in \mathbb N$, there exists a weak solution $(u_n,u_n')$ for (\ref{sys1})-(\ref{sys2}) on $[0,T]$ with initial data $(u_0^n,u_1^n)$.
We aim to show that $(u_n,u_n')$ converges to $(u,u')$ in the sense of (\ref{cdconverge}).

By (\ref{expm}) in Theorem \ref{thm1}, we have
\begin{align}   \label{pb-1}
&\frac{1}{2}\left(\|\nabla u_n(t)\|^2_2+\|u_n'(t)\|^2_2 \right)+\frac{1}{m+1}\norm{u_n(t)}_{m+1}^{m+1} +   \frac{1}{2} \int_0^t  \|u_n'(\tau)\|^{m+1}_{m+1} d\t   \notag\\\
&\leq  \left(    \frac{1}{2} (\|\nabla u_0^n\|_2^2 + \|u_1^n\|_2^2) + \frac{1}{m+1} \|u_0^n\|_{m+1}^{m+1}   + t\right)      e^{Ct} ,  \;\;   \text{for all}  \,\, t\in [0,T].
\end{align}
Also, one has
\begin{align}  \label{pb-2}
&\frac{1}{2}\left(\|\nabla u(t)\|^2_2+\|u'(t)\|^2_2 \right)+\frac{1}{m+1}\norm{u(t)}_{m+1}^{m+1} +   \frac{1}{2} \int_0^t  \|u'(\tau)\|^{m+1}_{m+1} d\tau \notag\\
&\leq    \left(    \frac{1}{2} (\|\nabla u_0\|_2^2 + \|u_1\|_2^2) + \frac{1}{m+1} \|u_0\|_{m+1}^{m+1}   + t\right)      e^{Ct},     \;\;   \text{for all}  \,\, t\in [0,T].
\end{align}

Notice $\|u\|_{H^1}^2 = \|\nabla u\|_2^2 + \|u\|_2^2  \leq   \|\nabla u\|_2^2 + C\|u\|_{m+1}^{m+1} + 1$ for $m\geq 1$. Then, 
since $\lim_{n \rightarrow \infty}\|u_0^n - u_0\|_{H^1}=0$, $\lim_{n \rightarrow \infty}\|u_0^n - u_0\|_{m+1}=0$ and $\lim_{n \rightarrow \infty}\|u_1^n - u_1\|_2=0$, 
and on account of (\ref{pb-1})-(\ref{pb-2}), there exists $K>0$ such that
\begin{align}   \label{boundK}
&\|u_n(t)\|_{H^1}^2 + \|u_n'(t)\|^2_2 +    \norm{u_n(t)}_{m+1}^{m+1}  +     \int_0^t  \|u_n'(\tau)\|^{m+1}_{m+1} d\tau \notag\\
&+ \|u(t)\|_{H^1}^2 + \|u'(t)\|^2_2 +    \norm{u(t)}_{m+1}^{m+1}  +     \int_0^t  \|u'(\tau)\|^{m+1}_{m+1} d\tau  \leq K,
\end{align}
for all $t\in [0,T]$, for all $n\in \mathbb N$.

Denote $y_n = u_n  - u$. Then, for all $t\in [0,T]$,
\begin{align}   \label{unique-ene}
&\frac{1}{2}  \left(       \|\nabla y_n(t)\|_2^2 +    \| y_n'(t)\|_2^2       \right) +\int_0^t \int_{\mathbb T^3} ( |u_n'|^{m-1}u_n' - |u'|^{m-1}u'  ) y_n' dx d\tau   \notag\\
&= \frac{1}{2}  \left(        \|\nabla y_n(0)\|_2^2 +       \| y_n'(0)\|_2^2      \right) 
+\int_0^t \int_{\mathbb T^3} (   |u_n|^{p-1} u_n - |u|^{p-1} u)  y_n' dx d\tau.
\end{align}
In fact, (\ref{unique-ene}) can be established rigorously by employing the regularization procedure used in the proof of the energy identity in subsection \ref{sec-EI}.

Since 
\begin{align*}
\frac{1}{m+1}\left(\|y_n(t)\|_{m+1}^{m+1}  -  \|y_n(0)\|_{m+1}^{m+1} \right)=  \int_0^t \int_{\mathbb T^3} |y_n(\tau)|^{m-1} y_n(\tau) y_n'(\tau) dx d\tau, 
\end{align*}
then (\ref{unique-ene}) can be written as
\begin{align}   \label{unique-e1}
&\frac{1}{2}  \left(       \|\nabla y_n(t)\|_2^2 +    \| y_n'(t)\|_2^2       \right) +  \frac{1}{m+1} \|y_n(t)\|_{m+1}^{m+1}
 +\int_0^t \int_{\mathbb T^3} ( |u_n'|^{m-1}u_n' - |u'|^{m-1}u'  ) y_n' dx d\tau   \notag\\
&= \frac{1}{2}  \left(        \|\nabla y_n(0)\|_2^2 +       \| y_n'(0)\|_2^2      \right)    +  \frac{1}{m+1}  \|y_n(0)\|_{m+1}^{m+1}    \notag\\
&\hspace{0.2  in} +\int_0^t \int_{\mathbb T^3} (   |u_n|^{p-1} u_n - |u|^{p-1} u)  y_n' dx d\tau   +  \int_0^t \int_{\mathbb T^3} |y_n|^{m-1} y_n y_n' dx d\tau.
\end{align}

Note, there exists a constant $c_0>0$ such that $(|a|^{m-1}a -  |b|^{m-1}b)(a-b)\geq c_0 |a-b|^{m+1}$ for all $a,b\in \mathbb R$. Then, since $y_n=u_n - u$, we have
\begin{align}   \label{unique-e2}
 \int_0^t \int_{\mathbb T^3} ( |u_n'|^{m-1}u_n' - |u'|^{m-1}u'  ) y_n' dx d\tau   \geq     c_0 \int_0^t \int_{\mathbb T^3} |y_n'|^{m+1} dx d\tau.
 \end{align}
Also, by the H\"older's inequality and Young's inequality, one has
\begin{align} \label{unique-e3}
\int_0^t \int_{\mathbb T^3} |y_n|^{m-1} y_n y_n' dx d\tau   
 \leq   c_0 \int_0^t \int_{\mathbb T^3} |y'_n|^{m+1}  dx d\tau+ C \int_0^t \int_{\mathbb T^3} |y_n|^{m+1}  dx d\tau.
\end{align}

Applying (\ref{unique-e2})-(\ref{unique-e3}) to equality (\ref{unique-e1}) yields
\begin{align}   \label{unique-enef}
&\frac{1}{2}  \left(       \|\nabla y_n(t)\|_2^2 +    \| y_n'(t)\|_2^2       \right) +  \frac{1}{m+1} \|y_n(t)\|_{m+1}^{m+1}   \notag\\
&= \frac{1}{2}  \left(        \|\nabla y_n(0)\|_2^2 +       \| y_n'(0)\|_2^2      \right)    +  \frac{1}{m+1}  \|y_n(0)\|_{m+1}^{m+1}    \notag\\
&\hspace{0.2  in} +\int_0^t \int_{\mathbb T^3} (   |u_n|^{p-1} u_n - |u|^{p-1} u)  y_n' dx d\tau   + C \int_0^t \|y_n(\tau)\|_{m+1}^{m+1} d\tau.
\end{align}

In the following, we estimate the integral $\int_0^t \int_{\mathbb T^3} (   |u_n|^{p-1} u_n - |u|^{p-1} u)  y_n' dx d\tau$.

For $1< p \leq 3$, by using H\"older's inequality and the imbedding $H^1\hookrightarrow L^6$, we have
\begin{align}   \label{sub}
&\int_0^t \int_{\mathbb T^3} (   |u_n|^{p-1} u_n - |u|^{p-1} u)  y_n' dx d\tau  \leq C\int_0^t \int_{\mathbb T^3} |y_n|   (|u_n|^{p-1} + |u|^{p-1})    |y_n'| dx d\tau  \notag\\
&\leq  C\int_0^t  \|y_n\|_6   \left(\|u_n\|_{3(p-1)}^{p-1} + \|u\|_{3(p-1)}^{p-1}   \right)   \|y_n'\|_2  d\tau  \leq C\int_0^t  \|y_n\|_6   (\|u_n\|_6^{p-1} + \|u\|_6^{p-1})   \|y_n'\|_2  d\tau \notag\\
&\leq C \int_0^t    (\|u_n\|_{H^1}^{p-1} + \|u\|_{H^1}^{p-1})  (\|y_n\|_{H^1}^2 + \|y_n'\|_2^2)    d\tau \leq C(K)  \int_0^t      (\|y_n\|_{H^1}^2 + \|y_n'\|_2^2) d\tau,
\end{align}
for all $t\in [0,T]$ and $n\in \mathbb N$, where the last inequality is owing to (\ref{boundK}).

Next we consider the ``supercritical" case $p>3$. Under such scenario, we apply integration by parts in time to convert $y_n'$ to $y_n$ in the integral $\int_0^t \int_{\mathbb T^3} (   |u_n|^{p-1} u_n - |u|^{p-1} u)  y_n' dx d\tau$. This idea originates from \cite{BL1} by Bociu and Lasiecka. 
Similar calculations have been performed in subsection \ref{AB} in the proof of existence of weak solutions. Indeed, $|u_n|^{p-1} u_n - |u|^{p-1} u = (u_n - u) \xi = y_n \xi$ where $|\xi| \leq C (|u_n|^{p-1}+ |u|^{p-1})$ and 
$|\xi'| \leq  C (|u_n|^{p-2}+ |u|^{p-2})  (|u_n'| + |u'|)  $.
Hence, by using integration by parts, we obtain
\begin{align}  \label{5terms}
&\int_0^t \int_{\mathbb T^3} (  |u_n|^{p-1} u_n - |u|^{p-1} u) y_n' dx d\tau     \notag\\
&=      \int_0^t \int_{\mathbb T^3} \xi  y_n y_n' dx d\tau  =  
\left[ \frac{1}{2}  \int_{\mathbb T^3}  \xi   y_n^2 dx \right]_0^t  -  \frac{1}{2} \int_0^t   \int_{\mathbb T^3} \xi' y_n^2  dx d\tau  \notag\\
&\leq    C \int_{\mathbb T^3}  \left(|u_n(t)|^{p-1} + |u(t)|^{p-1}\right) |y_n(t)|^2  dx  +    C \int_{\mathbb T^3}  \left(|u_n(0)|^{p-1} + |u(0)|^{p-1}  \right) |y_n(0)|^2  dx      \notag\\
& \hspace{0.5 in}  + C\int_0^t \int_{\mathbb T^3} \left(|u_n|^{p-2} + |u|^{p-2} \right) (|u_n'| + |u'|)y_n^2 dx d\tau,     \;\;  \text{for all} \,\,  t\in [0,T].\end{align}

We estimate each term on the right-hand side of (\ref{5terms}) as follows.

From Remark \ref{remark2}, we know that Case I and Case II in the assumption of Theorem \ref{thm2} can be combined as $p\leq \min\{ \frac{2}{3} m + \frac{5}{3}, m\}$. Thus $\frac{3(p-1)}{2}\leq m+1$. Then, using H\"older's inequality,
\begin{align}    \label{I-3}
&\int_{\mathbb T^3}(|u_n(0)|^{p-1}+|u(0)|^{p-1})|y_n(0)|^2 dx  
\leq \left(\norm{u_n(0)}^{p-1}_{\frac{3(p-1)}{2}}+\norm{u(0)}^{p-1}_{\frac{3(p-1)}{2}}\right)
\norm{y_n(0)}_{6}^2  \notag\\
&\leq C\left(\norm{u_n(0)}^{p-1}_{m+1}+\norm{u(0)}^{p-1}_{m+1}\right)
\norm{y_n(0)}_{H^1}^2
\leq   C(K) \norm{y_n(0)}_{H^1}^2,
\end{align}
due to (\ref{boundK}).

Furthermore, since $p\leq \frac{2}{3}m+\frac{5}{3}$, then $(p-2)\frac{3m+3}{2m-1}\leq m+1$. Therefore,
\begin{align}    \label{I-5}
&\int_0^t \int_{\mathbb T^3} (|u_n|^{p-2}+|u|^{p-2})(|u_n'|+|u'|)  y_n^2 dx d\t \notag\\
&\leq C\int_0^t \left(\norm{u_n}_{(p-2)\frac{3m+3}{2m-1}}^{p-2}+\norm{u}_{(p-2)\frac{3m+3}{2m-1}}^{p-2}\right) 
\left(\norm{u_n'}_{m+1}+\norm{u'}_{m+1}\right)    \norm{y_n}_{6}^2   d\tau    \notag\\
&\leq   C\int_0^t \left(\norm{u_n}_{m+1}^{p-2}+\norm{u}_{m+1}^{p-2}\right)
\left(\norm{u_n'}_{m+1}+\norm{u'}_{m+1}\right)    \norm{y_n}_{H^1}^2    d\tau \notag\\
&\leq C(K) \int_0^t \left(\norm{u_n'}_{m+1}+\norm{u'}_{m+1}\right)    \norm{y_n}_{H^1}^2  d\tau,  \;\;  \text{for all}  \,\,  t\in [0,T],
\end{align}
where we use (\ref{boundK}) to obtain the last inequality.

Finally, we consider $\int_{\mathbb T^3}(|u_n(t)|^{p-1}+|u(t)|^{p-1})|y_n(t)|^2 dx$.
We write out the estimate for $\int_{\mathbb T^3} |u_n(t)|^{p-1}   |y_n(t)|^2 dx$ only. 
The estimate for  $\int_{\mathbb T^3} |u(t)|^{p-1}   |y_n(t)|^2 dx$ is similar. Notice
\begin{align}   \label{split} 
&\int_{\mathbb T^3} |u_n(t)|^{p-1}   |y_n(t)|^2 dx    \leq C    \int_{\mathbb T^3} |u_n(t) - u_n(0)|^{p-1}   |y_n(t)|^2 dx    \notag\\
& \hspace{0.4 in} +   C    \int_{\mathbb T^3} |u_n(0)-u_0|^{p-1}  |y_n(t)|^2 dx 
+   C    \int_{\mathbb T^3} |u_0|^{p-1}  |y_n(t)|^2 dx.
\end{align}

Since $p\leq \frac{2}{3}m+\frac{5}{3}$, then $\frac{3(p-1)}{2(m+1)}\leq 1$. Therefore,
\begin{align}      \label{split1}
&\int_{\mathbb T^3} |u_n(t) - u_n(0)|^{p-1}   |y_n(t)|^2 dx   =   \int_{\mathbb T^3} \left|  \int_0^t u_n'(\tau)  d\tau \right|^{p-1}   |y_n(t)|^2 dx  \notag\\
 &\leq    \left(\int_{\mathbb T^3}       \left|  \int_0^t u_n'(\tau)  d\tau \right|^{\frac{3(p-1)}{2}} dx \right)^{2/3}    \|y_n(t)\|_6^2  \notag\\
 &\leq      C t^{\frac{m(p-1)}{m+1}} \left(\int_{\mathbb T^3}       \left(  \int_0^t |u_n'(\tau)|^{m+1}  d\tau \right)^{\frac{3(p-1)}{2(m+1)}} dx \right)^{2/3}  \|y_n(t)\|_{H^1}^2 \notag\\
 &\leq    C t^{\frac{m(p-1)}{m+1}} \left(\int_{\mathbb T^3}       \int_0^t |u_n'(\tau)|^{m+1}  d\tau dx \right)^{2/3}  \|y_n(t)\|_{H^1}^2 
 \leq C(K)    t^{\frac{m(p-1)}{m+1}}    \|y_n(t)\|_{H^1}^2,
   \end{align}
for all $t\in [0,T]$, due to (\ref{boundK}).

Next, we consider the integral $\int_{\mathbb T^3} |u_n(0)-u_0|^{p-1}  |y_n(t)|^2 dx$. 
Recall $u_n(0)=u_0^n  \rightarrow u_0$ in $L^{m+1}(\mathbb T^3)$. Note, the assumption that $p\leq \frac{2}{3} m + \frac{5}{3}$ implies $\frac{3(p-1)}{2}\leq m+1$. Then, by H\"older's inequality, one has
\begin{align}  \label{split2}
&\int_{\mathbb T^3} |u_n(0)-u_0|^{p-1}  |y_n(t)|^2 dx    =   \int_{\mathbb T^3} |u_0^n-u_0|^{p-1}  |y_n(t)|^2 dx
 \leq  \|u_0^n-u_0\|^{p-1}_{\frac{3(p-1)}{2}}     \|y_n(t)\|_6^2    \notag\\
&\leq   C\|u_0^n-u_0\|^{p-1}_{m+1}     \|y_n(t)\|_{H^1}^2   \leq  \epsilon  \|y_n(t)\|_{H^1}^2,   \;\; \text{for all}  \,\,  t\in [0,T],
\end{align}
for $n$ sufficiently large.

It remains to estimate the integral $  \int_{\mathbb T^3} |u_0|^{p-1}  |y_n(t)|^2 dx $.  
We notice
\begin{align}    \label{I-1}
\|y_n(t)\|_2^2 = \int_{\mathbb T^3} |y_n(t)|^2 dx
&=\int_{\mathbb T^3}\left(\Big|y_n(0)+\int_0^t y_n'(\t)d\tau \Big|^2 \right) dx   \notag\\
&\leq  \left(\norm{y_n(0)}_{2}^2 + t \int_0^t \norm{y_n'(\t)}^{2}_{2} d\t\right),       \;\;  \text{for all} \,\,  t\in [0,T].
\end{align}
Recall $\frac{3(p-1)}{2}\leq m+1$. Then, since $u_0 \in L^{m+1}(\mathbb T^3)  \subset   L^{\frac{3(p-1)}{2}}(\mathbb T^3) $,
there exists a smooth periodic function $\varphi$ such that $\|u_0-\varphi\|_{\frac{3(p-1)}{2}}^{p-1} \leq \epsilon$. Furthermore, since $\varphi$ is smooth on $\mathbb T^3$, there exists $C_{\epsilon}>0$ with $|\varphi(x)|  \leq C_{\epsilon}$ for all $x\in \mathbb T^3$. It follows that
\begin{align}   \label{split3}
&\int_{\mathbb T^3} |u_0|^{p-1}  |y_n(t)|^2 dx  \leq    C \int_{\mathbb T^3} |u_0 -\varphi |^{p-1}  |y_n(t)|^2 dx +  C\int_{\mathbb T^3} |\varphi|^{p-1}  |y_n(t)|^2 dx   \notag\\
&\leq C     \|u_0-\varphi\|_{\frac{3(p-1)}{2}}^{p-1}  \|y_n(t)\|_6^2  + C_{\epsilon}  \|y_n(t)\|_2^2      \notag\\
&\leq C \epsilon   \|y_n(t)\|_{H^1}^2   +  C_{\epsilon}    \left(\norm{y_n(0)}_{2}^2 + t \int_0^t \norm{y_n'(\t)}^{2}_{2} d\t\right),       \;\;  \text{for all} \,\,  t\in [0,T],
\end{align}
by virtue of (\ref{I-1}).

By substituting (\ref{split1}), (\ref{split2}) and (\ref{split3}) into (\ref{split}), we obtain, for $n$ sufficiently large, 
\begin{align}    \label{I-2}
&\int_{\mathbb T^3} |u_n(t)|^{p-1}   |y_n(t)|^2 dx    \notag\\
&\leq    \left(  C(K)  t^{\frac{m(p-1)}{m+1}}+C\epsilon \right)  \|y_n(t)\|_{H^1}^2
+C_{\epsilon}\left(\norm{y_n(0)}_{2}^2 + t\int_0^t \norm{y_n'(\t)}^{2}_{2} d\t\right),
\end{align}
for all $t\in [0,T]$. Using a similar calculation, $\int_{\mathbb T^3} |u(t)|^{p-1}   |y_n(t)|^2 dx$ has the same bound as (\ref{I-2}).

Now, applying the estimates (\ref{I-3}), (\ref{I-5}) and (\ref{I-2}) to (\ref{5terms}), it follows that
\begin{align}  \label{est-diff-f}
&\int_0^t \int_{\mathbb T^3} (  |u_n|^{p-1} u_n - |u|^{p-1} u) y_n' dx d\tau 
\leq  C(K,\epsilon)  \|y_n(0)\|_{H^1}^2   +  \left(  C(K)  t^{\frac{m(p-1)}{m+1}}+C\epsilon \right) \|y_n(t)\|_{H^1}^2  \notag\\
&  \hspace{0.3 in} +C(K)\int_0^t \left( \norm{u_n'(\tau)}_{m+1}+\norm{u'(\tau)}_{m+1} \right) \|y_n(\tau)\|_{H^1}^2 d\tau + C_{\epsilon} t\int_0^t \|y_n'(\tau)\|_2^2 d\tau,
\end{align}
for all $t\in [0,T]$ and for sufficiently large $n$, provided $p>3$.

We have finished estimating the integral $\int_0^t \int_{\mathbb T^3} (  |u_n|^{p-1} u_n - |u|^{p-1} u) y_n' dx d\tau$ for both of the subcritical case ($1\leq p\leq 3$) and the supercritical case ($p>3$). 
Then, since $\|y_n\|_{H^1}^2 =  \|\nabla y_n\|_2^2 +  \|y_n\|_2^2$, and using estimate (\ref{sub}),  (\ref{I-1}) and (\ref{est-diff-f}), we obtain from (\ref{unique-enef}) that
\begin{align}  \label{again}
&\frac{1}{2}  \left(       \|y_n(t)\|_{H^1}^2 +    \| y_n'(t)\|_2^2       \right) +\frac{1}{m+1} \|y_n(t)\|_{m+1}^{m+1}      \notag\\
&\leq   C(K,\epsilon)  \|y_n(0)\|_{H^1}^2   +   \frac{1}{2}\|y_n'(0)\|_2^2 +  \frac{1}{m+1} \|y_n(0)\|_{m+1}^{m+1}  +  \left(  C(K)  t^{\frac{m(p-1)}{m+1}}+C\epsilon \right) \|y_n(t)\|_{H^1}^2  \notag\\
& \hspace{0.3 in}+C(K)\int_0^t \left( \norm{u_n'(\tau)}_{m+1}+\norm{u'(\tau)}_{m+1}  +1  \right) \|y_n(\tau)\|_{H^1}^2 d\tau    \notag\\
&   \hspace{0.3 in} +   C(T,K,\epsilon)\int_0^t \|y_n'(\tau)\|_2^2 d\tau + C\int_0^t \|y_n(\tau)\|_{m+1}^{m+1} d\tau,
\end{align}
for all $t\in [0,T]$.

Then, by choosing $T_0\in (0,T]$ and $\epsilon>0$ sufficiently small such that   $ C(K)  T_0^{\frac{m(p-1)}{m+1}}+C\epsilon \leq \frac{1}{4}$, we obtain from (\ref{again}) that, for all $t\in [0,T_0]$,
\begin{align}  \label{again1}
&\frac{1}{4}  \left(       \|y_n(t)\|_{H^1}^2 +    \| y_n'(t)\|_2^2       \right) +\frac{1}{m+1} \|y_n(t)\|_{m+1}^{m+1}      \notag\\
&\leq   C(K,\epsilon)  \|y_n(0)\|_{H^1}^2   +   \frac{1}{2}\|y_n'(0)\|_2^2 +  \frac{1}{m+1} \|y_n(0)\|_{m+1}^{m+1}   \notag\\
& \hspace{0.3 in}+C(K)\int_0^t \left( \norm{u_n'(\tau)}_{m+1}+\norm{u'(\tau)}_{m+1}   +1    \right) \|y_n(\tau)\|_{H^1}^2 d\tau    \notag\\
&   \hspace{0.3 in} +   C(T,K,\epsilon)\int_0^t \|y_n'(\tau)\|_2^2 d\tau + C\int_0^t \|y_n(\tau)\|_{m+1}^{m+1} d\tau.
\end{align}

By virtue of (\ref{boundK}), we can apply the Gr\"onwall's inequality to (\ref{again1}) to conclude
\begin{align}  \label{unique-last}
&\|y_n(t)\|_{H^1}^2 +    \| y_n'(t)\|_2^2  + \|y_n(t)\|_{m+1}^{m+1}    \notag\\
&\leq C(K, T, \epsilon)   \left(\|y_n(0)\|_{H^1}^2   +   \|y_n'(0)\|_2^2   +  \|y_n(0)\|_{m+1}^{m+1}   \right), \;\;\text{for all}\,\, t\in [0,T_0].
\end{align}
Since $$\lim_{n\rightarrow \infty} ( \|y_n(0)\|_{H^1}^2 +  \|y_n'(0)\|_2^2  + \|y_n(0)\|_{m+1}^{m+1} )
 =   \lim_{n\rightarrow \infty} ( \|u_0^n - u_0\|_{H^1}^2 +  \|u_1^n-u_1\|_2^2  +    \|u_0^n - u_0\|_{m+1}^{m+1} )     =0 ,$$
then (\ref{unique-last}) implies that
$\lim_{n\rightarrow \infty} \left[ \sup_{t\in [0,T_0]} \left( \|y_n(t)\|_{H^1}^2 +   \| y_n'(t)\|_2^2  + \|y_n(t)\|_{m+1}^{m+1}   \right)\right] =0.$
By iterating the above procedure for finitely many times, we obtain
$$\lim_{n\rightarrow \infty} \left[ \sup_{t\in [0,T]} \left( \|y_n(t)\|_{H^1}^2 +   \| y_n'(t)\|_2^2  + \|y_n(t)\|_{m+1}^{m+1}   \right)\right] =0.$$
This completes the proof for the continuous dependence on initial data as well as the uniqueness of weak solutions.

\vspace{0.2 in}

\noindent {\bf Acknowledgment.} 
The author wishes to thank Mohammad Rammaha, Edriss Titi and Daniel Toundykov for helpful discussions.
\bibliographystyle{amsplain}

\end{document}